# FROM $\varepsilon$-ENTROPY TO KL-ENTROPY: ANALYSIS OF MINIMUM INFORMATION COMPLEXITY DENSITY ESTIMATION

By Tong Zhang

*Yahoo Research*

We consider an extension of $\varepsilon$-entropy to a KL-divergence based complexity measure for randomized density estimation methods. Based on this extension, we develop a general information-theoretical inequality that measures the statistical complexity of some deterministic and randomized density estimators. Consequences of the new inequality will be presented. In particular, we show that this technique can lead to improvements of some classical results concerning the convergence of minimum description length and Bayesian posterior distributions. Moreover, we are able to derive clean finite-sample convergence bounds that are not obtainable using previous approaches.

**1. Introduction.** The purpose of this paper is to study a class of complexity minimization based density estimation methods using a generalization of $\varepsilon$-entropy, which has become a central technical tool in the traditional finite-sample convergence analysis. Specifically, we derive a simple yet general information-theoretical inequality that can be used to measure the convergence of this very basic inequality.

We shall first introduce basic notation used in the paper. Consider a sample space $\mathcal{X}$ and a measure $\mu$ on $\mathcal{X}$ (with respect to some $\sigma$-field). In statistical inference, nature picks a probability measure $Q$ on $\mathcal{X}$ which is unknown. We assume that $Q$ has a density $q$ with respect to $\mu$. In density estimation, we consider a set of probability densities $p(\cdot|\theta)$ (with respect to $\mu$ on $\mathcal{X}$) indexed by $\theta \in \Gamma$. Without causing any confusion, we may also occasionally denote the model family $\{p(\cdot|\theta) : \theta \in \Gamma\}$ by the same symbol $\Gamma$. Throughout this paper, we always denote the true underlying density by $q$, and we do not assume that $q$ belongs to the model class $\Gamma$. Given $\Gamma$, our goal is to select a density $p(\cdot|\theta) \in \Gamma$ based on the observed data $X = \{X_1, \ldots, X_n\} \in \mathcal{X}^n$, such that $p(\cdot|\theta)$ is as close to $q$ as possible when measured by a certain distance function (which we shall specify later).

In the framework considered in this paper, we assume that there is a prior distribution $d\pi(\theta)$ on the parameter space $\Gamma$ that is independent of the observed data. For notational simplicity, we shall call any observation $X$ dependent probability density $\hat{w}_X(\theta)$ on $\Gamma$ (measurable on $\mathcal{X}^n \times \Gamma$) with respect to $d\pi(\theta)$ a *posterior randomization measure*, or simply a posterior. In particular, a posterior randomization measure in our sense is not limited to a *Bayesian posterior distribution*, which has a very specific meaning. We are interested in the density estimation performance of randomized estimators that draw $\theta$ according to posterior randomization measures $\hat{w}_X(\theta)$ obtained from a class of density estimation schemes. We should note that in this framework, our density estimator is completely characterized by the associated posterior $\hat{w}_X(\theta)$.

The paper is organized as follows. In Section 2, we introduce a generalization of $\varepsilon$-entropy for randomized estimation methods, which we call *KL-entropy*. Then a fundamental information-theoretical inequality, which forms the basis of our approach, will be obtained. Section 3 introduces the general information complexity minimization (ICM) density estimation formulation, where we derive various finite-sample convergence bounds using the fundamental information-theoretical inequality established earlier. Sections 4 and 5 apply the analysis to the case of minimum description length (MDL) estimators and to the convergence of Bayesian posterior distributions. In particular, we are able to simplify and improve most results in [1] as well as various recent analysis on the consistency and concentration of Bayesian posterior distributions. Some concluding remarks will be presented in Section 6.

Throughout this paper, we ignore the measurability issue, and assume that all quantities appearing in the derivations are measurable. Similarly to empirical process theory [14], the analysis can also be written in the language of outer-expectations, so that the measurability requirement imposed in this paper can be relaxed.

**2. The basic information-theoretical inequality.** In this section we introduce an information-theoretical complexity measure of randomized estimators represented as posterior randomization measures. As we shall see, this quantity directly generalizes the concept of $\varepsilon$-entropy for deterministic estimators. We also develop a simple yet very general information-theoretical inequality, which bounds the convergence behavior of an arbitrary randomized estimator using the introduced complexity measure. This inequality is the foundation of the approach introduced in this paper.

DEFINITION 2.1. Consider a probability density $w(\cdot)$ on $\Gamma$ with respect to $\pi$. The KL-divergence $D_{\mathrm{KL}}(w\,d\pi||d\pi)$ is defined as

$$D_{\mathrm{KL}}(w\,d\pi||d\pi) = \int_\Gamma w(\theta)\ln w(\theta)\,d\pi(\theta).$$

For any posterior randomization measure $\hat{w}_X$, we define its KL-entropy with respect to $\pi$ as $D_{\mathrm{KL}}(\hat{w}_X\,d\pi||d\pi)$.

Note that $D_{\mathrm{KL}}(w\,d\pi||d\pi)$ may not always be finite. However, it is always nonnegative.

KL-divergence is a rather standard information-theoretical concept. In this section we show that it can be used to measure the complexity of a randomized estimator. We can immediately see that the quantity directly generalizes the concept of $\varepsilon$-entropy on an $\varepsilon$-net; assuming that we have $N$ points in an $\varepsilon$-net, we may consider a prior that puts a mass of $1/N$ on every point. It is easy to see that any deterministic estimator in the $\varepsilon$-net can be regarded as a randomized estimator that is concentrated on one of the $N$ points with posterior weight $N$ (and weight of zero elsewhere). Clearly this estimator has a KL-entropy of $\ln N$, which is essentially the $\varepsilon$-entropy. In fact, it is also easy to verify that any randomized estimator on the $\varepsilon$-net has a KL-entropy bounded by its $\varepsilon$-entropy $\ln N$. Therefore $\varepsilon$-entropy is the worst-case KL-entropy on an $\varepsilon$-net with a uniform prior.

The concept of $\varepsilon$-entropy can be regarded as a notion to measure the complexity of an explicit discretization, usually for a deterministic estimator on a discrete $\varepsilon$-net. The concept of KL-entropy can be regarded as a notation to measure the complexity of a randomized estimation method, where the discretization is done implicitly through randomization with respect to an arbitrary prior. This difference is important for practical purposes since it is usually impossible (or very difficult) to perform computation on an explicitly discretized $\varepsilon$-net. Therefore estimators based on $\varepsilon$-nets are often of theoretical interest only. However, it is often feasible to draw samples from a posterior randomization measure with respect to a continuous prior by using standard Monte Carlo techniques. Therefore randomized estimation methods are potentially useful for practical problems.

Since KL-entropy allows nonuniform priors, the concept can directly characterize local adaptivity of randomized estimators when we put more prior mass in certain regions of the model family. In contrast, $\varepsilon$-entropy is a notation that tries to treat every part of the space equally, which may not give the best possible results. For example, for convergence of posterior distributions, the fact that entropy conditions are not always the most appropriate was pointed out in [4], pages 522–523. The issue of adaptivity (and related nonuniform prior) cannot be directly addressed with $\varepsilon$-entropy. In the literature, one has to employ additional techniques such as peeling (e.g., see [13]) for this purpose. As a comparison, the ability to use a nonuniform prior directly in our analysis is conceptually useful. Putting a large prior mass in a certain region indicates that we want to achieve a more accurate estimate in that region, in exchange for slower convergence in a region with smaller prior mass. The prior structure reflects our belief that the true density is more likely to have a certain form than some alternative forms. Therefore the theoretical analysis should also imply a more accurate estimate when we are lucky enough to guess the true density $q$ correctly by putting a large prior mass around it. As we will see later, finite-sample convergence bounds derived in this paper using KL-entropy have this behavior.

Next we prove a simple information-theoretical inequality using the KL-entropy of randomized estimators, which forms the basis of our analysis. For a real-valued function $f(\theta)$ on $\Gamma$, we denote by $\mathbf{E}_\pi f(\theta)$ the expectation of $f(\cdot)$ with respect to $\pi$. Similarly, for a real-valued function $\ell(x)$ on $\mathcal{X}$, we denote by $\mathbf{E}_q \ell(x)$ the expectation of $\ell(\cdot)$ with respect to the true underlying distribution $q$. We also use $\mathbf{E}_X$ to denote the expectation with respect to the observation $X$ ($n$ independent samples from $q$).

The key ingredient of our analysis using KL-entropy is a well-known convex duality, which has already been used in some recent machine learning papers to study sample complexity bounds. For example, see [8, 11]. For completeness, we include a simple information-theoretical proof.

PROPOSITION 2.1.    *Assume that $f(\theta)$ is a measurable real-valued function on $\Gamma$, and $w(\theta)$ is a density with respect to $\pi$; we have*

$$\mathbf{E}_\pi w(\theta) f(\theta) \le D_{\mathrm{KL}}(w\,d\pi \| d\pi) + \ln \mathbf{E}_\pi \exp(f(\theta)).$$

PROOF.    We assume that $\mathbf{E}_\pi \exp(f(\theta)) < \infty$; otherwise the bound is trivial. Consider $v(\theta) = \exp(f(\theta))/\mathbf{E}_\pi \exp(f(\theta))$. Since $\mathbf{E}_\pi v(\theta) = 1$, we can regard it as a density with respect to $\pi$. Using this definition, it is easy to verify that the inequality in Proposition 2.1 can be rewritten equivalently as

$$\mathbf{E}_\pi w(\theta) \ln w(\theta) + \ln \mathbf{E}_\pi \exp(f(\theta)) - \mathbf{E}_\pi w(\theta) f(\theta) = D_{\mathrm{KL}}(w\,d\pi \| v\,d\pi) \ge 0,$$

which is a well-known information-theoretical inequality, and follows easily from Jensen's inequality.    □

The main technical result which forms the basis of the paper is given by the following lemma, where we assume that $\hat{w}_X(\theta)$ is a posterior (represented as a density with respect to $\pi$ that depends on $X$ and is measurable on $\mathcal{X}^n \times \Gamma$).

LEMMA 2.1. *Consider any posterior $\hat{w}_X(\theta)$. Let $\alpha$ and $\beta$ be two real numbers. The following inequality holds for all measurable real-valued functions $L_X(\theta)$ on $\mathcal{X}^n \times \Gamma$:*

$$\mathbf{E}_X \exp[\mathbf{E}_\pi \hat{w}_X(\theta)(L_X(\theta) - \alpha \ln \mathbf{E}_X e^{\beta L_X(\theta)}) - D_{\mathrm{KL}}(\hat{w}_X \, d\pi || d\pi)]$$
$$\leq \mathbf{E}_\pi \frac{\mathbf{E}_X e^{L_X(\theta)}}{\mathbf{E}_X^\alpha e^{\beta L_X(\theta)}},$$

*where $\mathbf{E}_X$ is the expectation with respect to the observation $X$.*

PROOF. From Proposition 2.1, we obtain

$$\hat{L}(X) = \mathbf{E}_\pi \hat{w}_X(\theta)(L_X(\theta) - \alpha \ln \mathbf{E}_X e^{\beta L_X(\theta)}) - D_{\mathrm{KL}}(\hat{w}_X \, d\pi || d\pi)$$
$$\leq \ln \mathbf{E}_\pi \exp(L_X(\theta) - \alpha \ln \mathbf{E}_X e^{\beta L_X(\theta)}).$$

Now applying Fubini's theorem to interchange the order of integration, we have

$$\mathbf{E}_X e^{\hat{L}(X)} \leq \mathbf{E}_X \mathbf{E}_\pi e^{L_X(\theta) - \alpha \ln \mathbf{E}_X \exp(\beta L_X(\theta))} = \mathbf{E}_\pi \frac{\mathbf{E}_X e^{L_X(\theta)}}{\mathbf{E}_X^\alpha e^{\beta L_X(\theta)}}. \qquad \square$$

REMARK 2.1. The importance of the above inequality is that the left-hand side is a quantity that involves an arbitrary posterior randomization measure $\hat{w}_X \, d\pi$. The right-hand side is a numerical constant independent of the estimator $\hat{w}_X$. Therefore the inequality gives a bound that can be applied to an arbitrary randomized estimator. The remaining issue is merely how to interpret the resulting bound, which we shall focus on later in this paper.

REMARK 2.2. The main technical ingredients of the proof are motivated from techniques in the recent machine learning literature. The general idea for analyzing randomized estimators using Fubini's theorem and decoupling was already in [17]. The specific decoupling mechanism using Proposition 2.1 appeared in [3]; see [8, 11] for related problems. A simplified form of Lemma 2.1 was used in [18] to analyze Bayesian posterior distributions.

The following bound is a straightforward consequence of Lemma 2.1. Note that for density estimation, the loss $\ell_\theta(x)$ has the form of $\ell(p(x|\theta))$, where $\ell(\cdot)$ is a scaled log-loss.

THEOREM 2.1. *We use the notation of Lemma 2.1. Let $X = \{X_1, \ldots, X_n\}$ be $n$-samples that are independently drawn from $q$. Consider a measurable function $\ell_\theta(x) : \Gamma \times \mathcal{X} \to R$, and real numbers $\alpha$ and $\beta$, and define*

$$c_n(\alpha, \beta) = \frac{1}{n} \ln \mathbf{E}_\pi \left( \frac{\mathbf{E}_q e^{-\ell_\theta(x)}}{\mathbf{E}_q^\alpha e^{-\beta \ell_\theta(x)}} \right)^n.$$

*Then $\forall t$, the following event holds with probability at least $1 - \exp(-t)$:*

$$-\alpha \mathbf{E}_\pi \hat{w}_X(\theta) \ln \mathbf{E}_q e^{-\beta \ell_\theta(x)}$$

$$\leq \frac{\mathbf{E}_\pi \hat{w}_X(\theta) \sum_{i=1}^n \ell_\theta(X_i) + D_{\mathrm{KL}}(\hat{w}_X \, d\pi || d\pi) + t}{n} + c_n(\alpha, \beta).$$

*Moreover, we have the expected risk bound*

$$-\alpha \mathbf{E}_X \mathbf{E}_\pi \hat{w}_X(\theta) \ln \mathbf{E}_q e^{-\beta \ell_\theta(x)}$$

$$\leq \mathbf{E}_X \frac{\mathbf{E}_\pi \hat{w}_X(\theta) \sum_{i=1}^n \ell_\theta(X_i) + D_{\mathrm{KL}}(\hat{w}_X \, d\pi || d\pi)}{n} + c_n(\alpha, \beta).$$

Proof.   We use the notation of Lemma 2.1, with $L_X(\theta) = -\sum_{i=1}^n \ell_\theta(X_i)$. If we define

$$\hat{L}(X) = \mathbf{E}_\pi \hat{w}_X(\theta)(L_X(\theta) - \alpha \ln \mathbf{E}_X e^{\beta L_X(\theta)}) - D_{\mathrm{KL}}(\hat{w}_X \, d\pi || d\pi)$$

$$= \mathbf{E}_\pi \hat{w}_X(\theta) \left( -\sum_{i=1}^n \ell_\theta(X_i) - n\alpha \ln \mathbf{E}_q e^{-\beta \ell_\theta(x)} \right) - D_{\mathrm{KL}}(\hat{w}_X \, d\pi || d\pi),$$

then by Lemma 2.1 we have $\mathbf{E}_X e^{\hat{L}(X)} \leq e^{nc_n(\alpha,\beta)}$. This implies $\forall \varepsilon$: $e^\varepsilon P(\hat{L}(X) > \varepsilon) \leq e^{nc_n(\alpha,\beta)}$. Now given any $t$, and letting $\varepsilon = t + nc_n(\alpha, \beta)$, we obtain

$$e^{t + nc_n(\alpha,\beta)} P(\hat{L}(X) > t + nc_n(\alpha,\beta)) \leq e^{nc_n(\alpha,\beta)}.$$

That is, with probability at least $1 - e^{-t}$, $\hat{L}(X) \leq nc_n(\alpha, \beta) + t$. By rearranging the equation, we establish the first inequality of the theorem.

To prove the second inequality, we still start with $\mathbf{E}_X e^{\hat{L}(X)} \leq e^{nc_n(\alpha,\beta)}$ from Lemma 2.1. From Jensen's inequality with the convex function $e^x$, we obtain $e^{\mathbf{E}_X \hat{L}(X)} \leq \mathbf{E}_X e^{\hat{L}(X)} \leq e^{nc_n(\alpha,\beta)}$. That is, $\mathbf{E}_X \hat{L}(X) \leq nc(\alpha, \beta)$. By rearranging the equation, we obtain the desired bound.   □

Remark 2.3.   The special case of Theorem 2.1 with $\alpha = \beta = 1$ is very useful since in this case the term $c_n(\alpha, \beta)$ vanishes. In fact, in order to obtain the correct rate of convergence for nonparametric problems, it is sufficient to choose $\alpha = \beta = 1$. The more complicated case with general $\alpha$ and $\beta$ is only needed for parametric problems, where we would like to obtain a convergence rate of the order $O(1/n)$. In such cases the choice of $\alpha = \beta = 1$ would lead to a rate of $O(\ln n/n)$, which is suboptimal.

**3. Information complexity minimization.** Let $S$ be a predefined set of densities on $\Gamma$ with respect to the prior $\pi$. We consider a general information complexity minimization estimator,

$$(1) \qquad \hat{w}_X^S = \underset{w \in S}{\arg\min} \left[ -\mathbf{E}_\pi w(\theta) \sum_{i=1}^n \ln p(X_i|\theta) + \lambda D_{\mathrm{KL}}(w \, d\pi || d\pi) \right].$$

Given the true density $q$, if we define

$$(2) \qquad \hat{R}_\lambda(w) = \frac{1}{n} \mathbf{E}_\pi w(\theta) \sum_{i=1}^n \ln \frac{q(X_i)}{p(X_i|\theta)} + \frac{\lambda}{n} D_{\mathrm{KL}}(w\, d\pi||d\pi),$$

then it is clear that

$$\hat{w}_X^S = \underset{w \in S}{\arg\min}\, \hat{R}_\lambda(w).$$

The above estimation procedure finds a randomized estimator by minimizing the regularized empirical risk $\hat{R}_\lambda(w)$ among all possible densities with respect to the prior $\pi$ in a predefined set $S$. The purpose of this section is to study the performance of this estimator using Theorem 2.1. For simplicity, we shall only study the expected performance using the second inequality, although similar results can be obtained using the first inequality (which leads to exponential probability bounds).

One may define the true risk of $w$ by replacing the empirical expectation in (2) with the true expectation with respect to $q$:

$$(3) \qquad R_\lambda(w) = \mathbf{E}_\pi w(\theta) D_{\mathrm{KL}}(q||p(\cdot|\theta)) + \frac{\lambda}{n} D_{\mathrm{KL}}(w\, d\pi||d\pi),$$

where $D_{\mathrm{KL}}(q||p) = \mathbf{E}_q \ln(q(x)/p(x))$ is the KL-divergence between $q$ and $p$. The information complexity minimizer in (1) can be regarded as an approximate solution to (3) using empirical expectation.

Using empirical process techniques, one can typically expect to bound $R_\lambda(w)$ in terms of $\hat{R}_\lambda(w)$. Unfortunately, it does not work in our case since $D_{\mathrm{KL}}(q||p)$ is not well defined for all $p$. This implies that as long as $w$ has nonzero concentration around a density $p$ with $D_{\mathrm{KL}}(q||p) = +\infty$, then $R_\lambda(w) = +\infty$. Therefore we may have $R_\lambda(\hat{w}_X^S) = +\infty$ with nonzero probability even when the sample size approaches infinity.

A remedy is to use a distance function that is always well defined. In statistics, one often considers the $\rho$-divergence for $\rho \in (0,1)$, which is defined as

$$(4) \qquad D_\rho(q||p) = \frac{1}{\rho(1-\rho)} \mathbf{E}_q \left[ 1 - \left(\frac{p(x)}{q(x)}\right)^\rho \right].$$

This divergence is always well defined and $D_{\mathrm{KL}}(q||p) = \lim_{\rho \to 0} D_\rho(q||p)$. In the statistical literature, convergence results were often specified under the squared Hellinger distance ($\rho = 0.5$). In this paper we specify convergence results with general $\rho$. We shall mention that bounds derived in this paper will become trivial when $\rho \to 0$. This is consistent with the above discussion since $R_\lambda$ (corresponding to $\rho = 0$) may not converge at all. However, under additional assumptions, such as the boundedness of $q/p$, $D_{\mathrm{KL}}(q||p)$ exists and can be bounded using the $\rho$-divergence $D_\rho(q||p)$.

A concept related to the $\rho$-divergence in (4) is the Rényi entropy introduced in [9]. The notion has been widely used in information theory. Up to a scaling factor, it can be defined as

$$D_\rho^{\mathrm{Re}}(q||p) = -\frac{1}{\rho(1-\rho)} \ln \mathbf{E}_q \left( \frac{p(x)}{q(x)} \right)^\rho.$$

Note that the standard definition of Rényi entropy in the literature is $\rho D_\rho^{\mathrm{Re}}(q||p)$. We employ a scaled version in this paper for compatibility with our $\rho$-divergence definition. Using the inequality $1-x \leq -\ln x \leq x^{-1}-1$ ($x \in [0,1]$), we can see that $\forall\, p, q$

$$D_\rho(q||p) \leq D_\rho^{\mathrm{Re}}(q||p) \leq \frac{D_\rho(q||p)}{1-\rho(1-\rho)D_\rho(q||p)}.$$

The following bounds imply that up to a constant, the $\rho$-divergence with any $\rho \in (0,1)$ is equivalent to the squared Hellinger distance. Therefore a convergence bound in any $\rho$-divergence implies a convergence bound of the same rate in the Hellinger distance.

PROPOSITION 3.1.  *We have the following inequalities $\forall\, \rho \in [0,1]$:*

$$\max(\rho, 1-\rho)D_\rho(q||p) \geq \tfrac{1}{2}D_{1/2}(q||p) \geq \min(\rho, 1-\rho)D_\rho(q||p).$$

PROOF.  We prove the first half of the two inequalities. Due to the symmetry $D_\rho(q||p) = D_{1-\rho}(p||q)$, we only need to consider the case $\rho \leq 1/2$. The proof of the second half (with $\rho \geq 1/2$) is identical except that the sign in the Taylor expansion step is reversed.

We use Taylor expansion. Let $x = \frac{p^{1/2}-q^{1/2}}{q^{1/2}}$; then $x \geq -1$, and there exists $\xi > -1$ such that

$$(1+x)^{2\rho} = 1 + 2\rho x + \rho(2\rho-1)(1+\xi)^{2\rho-2}x^2 \leq 1 + 2\rho x.$$

Now taking expectation with respect to $q$, we obtain

$$\mathbf{E}_q \left( \frac{p}{q} \right)^\rho = \mathbf{E}_q \left( 1 + \frac{p^{1/2}-q^{1/2}}{q^{1/2}} \right)^{2\rho} \leq 1 + 2\rho \mathbf{E}_q \frac{p^{1/2}-q^{1/2}}{q^{1/2}}.$$

By rearranging the equation, we obtain $2\rho(\tfrac{1}{4}D_{1/2}(q||p)) \leq \rho(1-\rho)D_\rho(q||p)$. □

3.1. *A general convergence bound.*  The following theorem is a consequence of Theorem 2.1. Most of our later discussion can be considered as interpretation of this theorem under different conditions.

THEOREM 3.1. *Consider the estimator $\hat{w}_X^S$ defined in* (1). *Let $\alpha > 0$. Then $\forall \rho \in (0,1)$ and $\gamma \geq \rho$ such that $\lambda' = \frac{\lambda\gamma - 1}{\gamma - \rho} \geq 0$, we have*

$$\mathbf{E}_X \mathbf{E}_\pi \hat{w}_X^S(\theta) D_\rho(q||p(\cdot|\theta)) \leq \mathbf{E}_X \mathbf{E}_\pi \hat{w}_X^S(\theta) D_\rho^{\mathrm{Re}}(q||p(\cdot|\theta))$$

$$\leq \frac{\gamma \inf_{w \in S} R_\lambda(w)}{\alpha \rho (1-\rho)} - \frac{\gamma - \rho}{\alpha \rho (1-\rho)} \mathbf{E}_X \hat{R}_{\lambda'}(\hat{w}_X^S)$$

$$+ \frac{c_{\rho,n}(\alpha)}{\alpha \rho (1-\rho)},$$

*where $c_{\rho,n}(\alpha) = \frac{1}{n} \ln \mathbf{E}_\pi \mathbf{E}_q^{(1-\alpha)n} (\frac{p(x|\theta)}{q(x)})^\rho = \frac{1}{n} \ln \mathbf{E}_\pi e^{-\rho(1-\rho)(1-\alpha)n D_\rho^{\mathrm{Re}}(q||p(\cdot|\theta))}$.*

PROOF. Consider an arbitrary data-independent density $w(\theta) \in S$ with respect to $\pi$. Using (4), we can obtain from Theorem 2.1 the chain of equations

$$\alpha \rho (1-\rho) \mathbf{E}_X \mathbf{E}_\pi \hat{w}_X^S(\theta) D_\rho(q||p(\cdot|\theta))$$

$$\leq \alpha \rho (1-\rho) \mathbf{E}_X \mathbf{E}_\pi \hat{w}_X^S(\theta) D_\rho^{\mathrm{Re}}(q||p(\cdot|\theta))$$

$$= -\alpha \mathbf{E}_X \mathbf{E}_\pi \hat{w}_X^S(\theta) \ln \mathbf{E}_q \exp\left(-\rho \ln \frac{q(x)}{p(x|\theta)}\right)$$

$$\leq \mathbf{E}_X \left[ \rho \mathbf{E}_\pi \hat{w}_X^S \sum_{i=1}^n \frac{1}{n} \ln \frac{q(X_i)}{p(X_i|\theta)} + \frac{D_{\mathrm{KL}}(\hat{w}_X^S \, d\pi || d\pi)}{n} \right] + c_{\rho,n}(\alpha)$$

$$= \mathbf{E}_X [\gamma \hat{R}_\lambda(\hat{w}_X^S) + (\rho - \gamma) \hat{R}_{\lambda'}(\hat{w}_X^S)] + c_{\rho,n}(\alpha)$$

$$\leq \mathbf{E}_X [\gamma \hat{R}_\lambda(w) + (\rho - \gamma) \hat{R}_{\lambda'}(\hat{w}_X^S)] + c_{\rho,n}(\alpha)$$

$$= \gamma R_\lambda(w) - (\gamma - \rho) \mathbf{E}_X \hat{R}_{\lambda'}(\hat{w}_X^S) + c_{\rho,n}(\alpha),$$

where $R_\lambda(w)$ is defined in (3). Note that the first inequality uses the fact $-\ln(1-x) \geq x$. The second inequality follows from Theorem 2.1 with the choice $\ell_\theta(x) = \rho \ln \frac{q(x)}{p(x|\theta)}$ and $\beta = 1$. The third inequality follows from the definition of $\hat{w}_X^S$ in (1). □

REMARK 3.1. If $\gamma = \rho$ in Theorem 3.1, then we also require $\lambda\gamma = 1$, and let $\lambda' = 0$.

Although the bound in Theorem 3.1 looks complicated, the most important part on the right-hand side is the first term. The second term is only needed to handle the situation $\lambda \leq 1$. The requirement that $\gamma \geq \rho$ is to ensure that the second term is nonpositive. Therefore in order to apply the theorem, we only need to estimate a lower bound of $\hat{R}_{\lambda'}(\hat{w}_X^S)$, which (as we shall see later) is much easier than obtaining an upper bound. The third term is mainly included to get the correct convergence rate of $O(1/n)$ for parametric problems, and can be ignored for nonparametric problems. The effect of this term is quite similar to using localized $\varepsilon$-entropy in the empirical process approach for analyzing the maximum-likelihood method; for example, see [13]. As a comparison, the KL-entropy in the first term corresponds to the global $\varepsilon$-entropy.

Note that one can easily obtain a simplified bound from Theorem 3.1 by choosing specific parameters so that both the second term and the third term vanish:

COROLLARY 3.1. *Consider the estimator $\hat{w}_X^S$ defined in (1). Assume that $\lambda > 1$ and let $\rho = 1/\lambda$. We have*

$$\mathbf{E}_X \mathbf{E}_\pi \hat{w}_X^S(\theta) D_\rho^{\mathrm{Re}}(q||p(\cdot|\theta)) \leq \frac{1}{1-\rho} \inf_{w \in S} R_\lambda(w).$$

PROOF. We simply let $\alpha = 1$ and $\gamma = \rho$ in Theorem 3.1. □

An important observation is that for $\lambda > 1$, the convergence rate is solely determined by the quantity $\inf_{w \in S} R_\lambda(w)$, which we shall refer to as the *model resolvability* associated with $S$.

3.2. *Some consequences of Theorem 3.1.* In order to apply Theorem 3.1, we need to bound the quantity $\mathbf{E}_X \hat{R}_{\lambda'}(\hat{w}_X^S)$ from below. Some of these results can be found in the Appendix, and by using these results, we are able to obtain some refined bounds from Theorem 3.1.

COROLLARY 3.2. *Consider the estimator $\hat{w}_X^S$ defined in (1). Assume that $\lambda > 1$; then $\forall \rho \in (0, 1/\lambda]$*

$$\mathbf{E}_X \mathbf{E}_\pi \hat{w}_X^S(\theta) D_\rho^{\mathrm{Re}}(q||p(\cdot|\theta)) \leq \frac{1}{\rho(\lambda-1)} \inf_{w \in S} R_\lambda(w).$$

PROOF. We simply let $\alpha = 1$ and $\gamma = (1-\rho)/(\lambda-1)$ in Theorem 3.1. Note that in this case, $\lambda' = 1$, and hence by Lemma A.1 in the Appendix, we have $\mathbf{E}_X \hat{R}_{\lambda'}(\hat{w}_X^S) \geq 0$. □

Note that Lemma A.1 is only applicable for $\lambda' \geq 1$. If $\lambda' \leq 1$, then we need a discretization device which generalizes the upper $\varepsilon$-covering number concept used in [2] for showing the consistency (or inconsistency) of Bayesian posterior distributions:

DEFINITION 3.1. The $\varepsilon$-upper bracketing number of $\Gamma$, denoted by $N_{\mathrm{ub}}(\Gamma, \varepsilon)$, is the minimum number of nonnegative functions $\{f_j\}$ on $\mathcal{X}$ with respect to $\mu$ such that $\mathbf{E}_q(f_j/q) \leq 1 + \varepsilon$, and $\forall \theta \in \Gamma$, $\exists j$ such that $p(x|\theta) \leq f_j(x)$ a.e. $[\mu]$.

The discretization device which we shall use in this paper is based on the following definition.

DEFINITION 3.2. Given a set $\Gamma' \subset \Gamma$, we define its upper-bracketing radius as

$$r_{\mathrm{ub}}(\Gamma') = \int \sup_{\theta \in \Gamma'} p(x|\theta) \, d\mu(x) - 1.$$

An $\varepsilon$-upper discretization of $\Gamma$ consists of a covering of $\Gamma$ by countably many measurable subsets $\{\Gamma_j\}$ such that $\bigcup_j \Gamma_j = \Gamma$ and $r_{\mathrm{ub}}(\Gamma_j) \leq \varepsilon$.

Using this concept, we may combine the estimate in Lemma A.2 in the Appendix with Theorem 3.1, and obtain the following simplified bound for $\lambda = 1$. Similar results can also be obtained for $\lambda < 1$.

COROLLARY 3.3. Consider the estimator defined in (1). Let $\lambda = 1$. Consider an arbitrary covering $\{\Gamma_j\}$ of $\Gamma$. $\forall \rho \in (0,1)$ and $\forall \gamma \geq 1$, we have

$$\mathbf{E}_X \mathbf{E}_\pi \hat{w}_X^{\mathrm{S}}(\theta) D_\rho^{\mathrm{Re}}(q||p(\cdot|\theta))$$
$$\leq \frac{\gamma \inf_{w \in S} R_\lambda(w)}{\rho(1-\rho)} + \frac{\gamma - \rho}{\rho(1-\rho)n} \ln \sum_j \pi(\Gamma_j)^{(\gamma-1)/(\gamma-\rho)} (1 + r_{\mathrm{ub}}(\Gamma_j))^n.$$

In particular, if $\{\Gamma_j^\varepsilon\}$ is an $\varepsilon$-upper discretization of $\Gamma$, then

$$\mathbf{E}_X \mathbf{E}_\pi \hat{w}_X^{\mathrm{S}}(\theta) D_\rho^{\mathrm{Re}}(q||p(\cdot|\theta))$$
$$\leq \frac{\gamma \inf_{w \in S} R_\lambda(w)}{\rho(1-\rho)} + \frac{\gamma - \rho}{\rho(1-\rho)} \left[ \frac{\ln \sum_j \pi(\Gamma_j^\varepsilon)^{(\gamma-1)/(\gamma-\rho)}}{n} + \ln(1+\varepsilon) \right].$$

PROOF. We let $\alpha = 1$ in Theorem 3.1 and apply Lemma A.2. $\square$

Note that the above results immediately imply the following bound using $\varepsilon$-upper entropy by letting $\gamma \to 1$ with a finite $\varepsilon$-upper bracketing cover of size $N_{\mathrm{ub}}(\Gamma, \varepsilon)$ as the discretization:

$$\mathbf{E}_X \mathbf{E}_\pi \hat{w}_X^{\mathrm{S}}(\theta) D_\rho^{\mathrm{Re}}(q||p(\cdot|\theta)) \leq \frac{\inf_{w \in S} R_\lambda(w)}{\rho(1-\rho)} + \frac{1}{\rho} \inf_{\varepsilon > 0} \left[ \frac{\ln N_{\mathrm{ub}}(\Gamma, \varepsilon)}{n} + \ln(1+\varepsilon) \right].$$

It is clear that Corollary 3.3 is significantly more general. We are able to deal with an infinite cover as long as the decay of the prior $\pi$ is fast enough on an $\varepsilon$-upper discretization so that $\sum_j \pi(\Gamma_j^\varepsilon)^{(\gamma-1)/(\gamma-\rho)} < +\infty$.

3.3. *Weak convergence bound.* The case of $\lambda = 1$ is related to a number of important estimation methods in statistical applications. However, for an arbitrary prior $\pi$ without any additional assumption such as the fast decay condition in Corollary 3.3, it is impossible to establish any convergence rate result in terms of Hellinger distance using the model resolvability quantity alone, as in the case of $\lambda > 1$ (Corollary 3.2). See Section 4.4 for an example demonstrating this claim. However, one can still obtain a weaker convergence result in this case.

THEOREM 3.2. *Consider the estimator $\hat{w}_X^S$ defined in (1) with $\lambda = 1$. Then $\forall f : \mathcal{X} \to [-1, 1]$, we have*

$$\mathbf{E}_X \left| \mathbf{E}_\pi \hat{w}_X^S(\theta) \mathbf{E}_{p(\cdot|\theta)} f(x) - \frac{1}{n} \sum_{i=1}^n f(X_i) \right| \le 2A_n + \sqrt{2A_n},$$

*where $\mathbf{E}_{p(\cdot|\theta)} f(x) = \int f(x) p(x|\theta) \, d\mu(x)$ is the expectation with respect to $p(\cdot|\theta)$ on $\mathcal{X}$, and $A_n = \inf_{w \in S} \mathbf{E}_\pi R_\lambda(w) + \frac{\ln 2}{n}$.*

PROOF. The first half of the proof, leading to (5), is an application of Theorem 2.1. The second half is very similar to the proof of Theorem 3.1.

Let $g_\varepsilon(x) = 1 - \varepsilon f(x)$, and $h_\varepsilon(\theta, x) = \frac{q(x)}{p(x|\theta) g_\varepsilon(x)}$, where $\varepsilon \in (-1, 1)$ is a parameter to be determined later. Note that $g_\varepsilon(x) > 0$.

We consider an extension of $\Gamma$ to $\Gamma' = \Gamma \times \{\pm 1\}$. Let $\sigma = \pm 1$, and $\theta' = (\theta, \sigma) \in \Gamma'$. We define a prior $\pi'$ on $\Gamma'$ such that $\pi'((\theta, \sigma)) = 0.5\pi(\theta)$. For a posterior $\hat{w}_X^S(\theta)$ on $\Gamma$, we consider for $u = \pm 1$ a posterior $\hat{w}_{u,X}^S(\theta, \sigma)$ on $\Gamma'$ such that $\hat{w}_{u,X}^S(\theta, \sigma) = 2\hat{w}_X^S(\theta)$ when $\sigma = u$, and $\hat{w}_{u,X}^S(\theta, \sigma) = 0$ otherwise. Let $\alpha = \beta = 1$, $\ell_{\theta,\sigma}(x) = \ln h_{\sigma\varepsilon}(\theta, X_i)$. For all $u(X) \in \{\pm 1\}$, we apply Theorem 2.1 to the posterior $\hat{w}_{u(X),X}^S$, and obtain

$$-\mathbf{E}_X \mathbf{E}_\pi \hat{w}_X(\theta) \ln \mathbf{E}_q e^{-\ln h_{u\varepsilon}(\theta, x)}$$
$$\le \mathbf{E}_X \frac{\mathbf{E}_\pi \hat{w}_X(\theta) \sum_{i=1}^n \ln h_{u\varepsilon}(\theta, X_i) + D_{\mathrm{KL}}(\hat{w}_X \, d\pi || d\pi) + \ln 2}{n}.$$

Note that $\mathbf{E}_q e^{-\ln h_{u\varepsilon}(\theta, x)} = \mathbf{E}_{p(\cdot|\theta)} g_\varepsilon(x)$. Therefore if we let

$$\Delta_\varepsilon(X) = \mathbf{E}_\pi \hat{w}_X^S(\theta) \left( \sum_{i=1}^n \ln g_\varepsilon(X_i) - n \ln \mathbf{E}_{p(\cdot|\theta)} g_\varepsilon(x) \right),$$

then

$$(5) \qquad \mathbf{E}_X \Delta_{u(X)\varepsilon}(X) \le n \mathbf{E}_X \hat{R}_\lambda(\hat{w}_X^S) + \ln 2 \le n \inf_{w \in S} R_\lambda(w) + \ln 2,$$

where the second inequality follows from the definition of $\hat{w}_X^S$ in (1). This inequality plays the same role as Theorem 2.1 in the proof of Theorem 3.1.

Consider $x \le y < 1$. We have the inequalities (which follow from Taylor expansion)

$$x \le -\ln(1-x) \le x + \frac{x^2}{2(1-y)^2}.$$

This implies $\ln g_\varepsilon(x) \ge -\varepsilon f(x) - \frac{\varepsilon^2}{2(1-|\varepsilon|)^2}$ and $-\ln \mathbf{E}_{p(\cdot|\theta)} g_\varepsilon(x) \ge \varepsilon \mathbf{E}_{p(\cdot|\theta)} f(x)$. Therefore

$$\Delta_\varepsilon(X) \ge \varepsilon \mathbf{E}_\pi \hat{w}_X^S(\theta)\left(-\sum_{i=1}^n f(X_i) + n\mathbf{E}_{p(\cdot|\theta)}f(x)\right) - \frac{n\varepsilon^2}{2(1-|\varepsilon|)^2}.$$

Substitute into (5); we have

$$\mathbf{E}_X \sup_{u \in \{\pm 1\}} \left(u\varepsilon \mathbf{E}_\pi \hat{w}_X^S(\theta)\left(-\sum_{i=1}^n f(X_i) + n\mathbf{E}_{p(\cdot|\theta)}f(x)\right)\right) - \frac{n\varepsilon^2}{2(1-|\varepsilon|)^2}$$
$$\le n \inf_{w \in S} \mathbf{E}_X R_\lambda(w) + \ln 2.$$

Therefore we have

$$\mathbf{E}_X \left|\mathbf{E}_\pi \hat{w}_X^S(\theta)\left(-\sum_{i=1}^n f(X_i) + n\mathbf{E}_{p(\cdot|\theta)}f(x)\right)\right| \le \frac{n|\varepsilon|}{2(1-|\varepsilon|)^2} + \frac{nA_n}{|\varepsilon|}.$$

Let $|\varepsilon| = \sqrt{2A_n}/(\sqrt{2A_n}+1)$ and we obtain the desired bound. $\square$

Note that for all $f \in [-1, 1]$ the empirical average $\frac{1}{n}\sum_{i=1}^n f(X_i)$ converges to $\mathbf{E}_q f(x)$,

$$\mathbf{E}_X \left|\frac{1}{n}\sum_{i=1}^n f(X_i) - \mathbf{E}_q f(x)\right| \le \mathbf{E}_X^{1/2}\left(\frac{1}{n}\sum_{i=1}^n f(X_i) - \mathbf{E}_q f(x)\right)^2 \le \frac{1}{\sqrt{n}}.$$

It follows from Theorem 3.2 that

$$\mathbf{E}_X |\mathbf{E}_\pi \hat{w}_X^S(\theta)\mathbf{E}_{p(\cdot|\theta)}f(x) - \mathbf{E}_q f(x)| \le 2A_n + \sqrt{2A_n} + n^{-1/2}.$$

This means that as long as $\lim_n A_n = 0$, for all bounded functions $f(x) \in [-1, 1]$, the posterior average $\mathbf{E}_\pi \hat{w}_X^S(\theta)\mathbf{E}_{p(\cdot|\theta)}f(x)$ converges to $\mathbf{E}_q f(x)$ in probability. Since Theorem 3.2 uses the same weak topology as that in the usual definition of weak convergence of measures, we can interpret this result to mean the posterior average $\mathbf{E}_\pi \hat{w}_X^S(\theta)p(\cdot|\theta)$ converges weakly to $q$ in probability. In particular, by letting $f(x)$ be an indicator function for an arbitrary set $B \subset \mathcal{X}$, we obtain the consistency of the probability estimate. That is, the probability of $B$ under the posterior mean $\mathbf{E}_\pi \hat{w}_X^S(\theta)p(\cdot|\theta)$ converges to the probability of $B$ under $q$ (when $\lim_n A_n = 0$).

**4. Two-part code MDL on discrete net.** The minimum description length (MDL) method has been widely used in practice [10]. The two-part code MDL we consider here is the same as that of Barron and Cover [1]. In fact, results in this section improve those of Barron and Cover [1]. The MDL method considered in [1] can be regarded as a special case of information complexity minimization. The model space $\Gamma$ is countable: $\theta \in \Gamma = \{1, 2, \ldots\}$. We denote the corresponding models $p(x|\theta = j)$ by $p_j(x)$. The prior $\pi$ has the form $\pi = \{\pi_1, \pi_2, \ldots\}$ such that $\sum_j \pi_j = 1$, where we assume that $\pi_j > 0$ for each $j$. A randomized algorithm can be represented as a nonnegative weight vector $w = [w_j]$ such that $\sum_j \pi_j w_j = 1$.

MDL gives a deterministic estimator, which corresponds to the set of weights concentrated on any one specific point $k$. That is, we can select $S$ in (1), where each weight $w$ in $S$ corresponds to an index $k \in \Gamma$ such that $w_k = 1/\pi_k$ and $w_j = 0$ when $j \neq k$. It is easy to check that $D_{\mathrm{KL}}(w \, d\pi || d\pi) = \ln(1/\pi_k)$. The corresponding algorithm can thus be described as finding a probability density $p_{\hat{k}}$ with $\hat{k}$ obtained by

$$(6) \qquad \hat{k} = \arg\min_k \left[ \sum_{i=1}^n \ln \frac{1}{p_k(X_i)} + \lambda \ln \frac{1}{\pi_k} \right],$$

where $\lambda \geq 1$ is a regularization parameter. The first term corresponds to the description of the data, and the second term corresponds to the description of the model. The choice $\lambda = 1$ can be interpreted as minimizing the total description length, which corresponds to the standard MDL. The choice $\lambda > 1$ corresponds to heavier penalty on the model description, which makes the estimation method more stable. This modified MDL method was considered in [1] and the authors obtained results on the asymptotic rate of convergence. However, no simple finite-sample bound was obtained. For the case of $\lambda = 1$, only weak consistency was shown. In the following, we shall improve these results using the analysis presented in Section 3.

4.1. *Modified MDL under global entropy condition.* Consider the case $\lambda > 1$ in (6). We can obtain the following theorem from Corollary 3.2.

THEOREM 4.1. *Consider the estimator $\hat{k}$ defined in* (6). *Assume that $\lambda > 1$. Then $\forall \rho \in (0, 1/\lambda]$*

$$\mathbf{E}_X D_\rho(q||p_{\hat{k}}) \leq \mathbf{E}_X D_\rho^{\mathrm{Re}}(q||p_{\hat{k}}) \leq \frac{1}{\rho(\lambda - 1)} \inf_k \left[ D_{\mathrm{KL}}(q||p_k) + \frac{\lambda}{n} \ln \frac{1}{\pi_k} \right].$$

The term $r_{\lambda,n}(q) = \inf_k [D_{\mathrm{KL}}(q||p_k) + \frac{\lambda}{n} \ln \frac{1}{\pi_k}]$ is referred to as *index of resolvability* in [1]. They showed (Theorem 4) that $D_{1/2}(q||p_{\hat{k}}) = O_p(r_{\lambda,n}(q))$ when $\lambda > 1$, which is a direct consequence of Theorem 4.1.

Theorem 4.1 generalizes a result by Andrew Barron and Jonathan Li, which gave a similar inequality but only for the case of $\lambda = 2$ and $\rho = 1/2$. The result can be found in [7], Theorem 5.5, page 78. In particular, consider $\Gamma$ such that $|\Gamma| = N$ with uniform prior $\pi_j = 1/N$; one obtains a bound for the maximum likelihood estimate over $\Gamma$ (take $\lambda = 2$ and $\rho = 1/2$ in Theorem 4.1),

$$(7) \qquad \mathbf{E}_X D_{1/2}(q||p_{\hat{k}}) \leq 2 \inf_k \left[ D_{\mathrm{KL}}(q||p_k) + \frac{2}{n} \ln \frac{1}{N} \right].$$

Examples of indexes of resolvability for various function classes can be found in [1], which we shall not repeat in this paper. In particular, it is known that for nonparametric problems, with appropriate discretization the rate resulting from (7) matches the minimax rate, such as those in [16].

4.2. *Local entropy analysis.* Although the bound based on the index of resolvability in Theorem 4.1 is quite useful for nonparametric problems, see [1], it does not handle the parametric case satisfactorily. To see this, we consider a one-dimensional parameter family indexed by $\theta \in [0, 1]$, and we discretize the family using a uniform discrete net of size $N+1$, $\theta_j = j/N$ ($j = 0, \ldots, N$). In the following, we assume that $q$ is taken from the parametric family, and for some fixed $\rho$, both $D_\rho^{\mathrm{Re}}(q||p_k)$ and $D_{\mathrm{KL}}(q||p_k)$ are of the order $(\theta - \theta_k)^2$. That is, we assume that there exist constants $c_1$ and $c_2$ where

$$(8) \qquad c_1(\theta - \theta_k)^2 \leq D_\rho^{\mathrm{Re}}(q||p_k), \qquad D_{\mathrm{KL}}(q||p_k) \leq c_2(\theta - \theta_k)^2.$$

We will thus have $\inf_k D_{\mathrm{KL}}(q||p_k) \leq c_2 N^{-2}$, and the bound in (7), which relies on the index of resolvability, becomes $\mathbf{E}_X D_{1/2}(q||p_{\hat{k}}) \leq O(N^{-2}) + \frac{4}{n} \ln \frac{1}{N+1}$. Now by choosing $N = O(n^{-1/2})$, we obtain a suboptimal convergence rate $\mathbf{E}_X D_{1/2}(q||p_{\hat{k}}) \leq O(\ln n/n)$. Note that convergence rates established in [1] for parametric examples are also of the order $O(\ln n/n)$.

The main reason for this suboptimality is that the complexity measure $O(\ln N)$ or $O(-\ln \pi_k)$ corresponds to the globally defined entropy. However, readers who are familiar with the empirical process theory know that the rate of convergence of the maximum-likelihood estimate is determined by local entropy mentioned in [5]. For nonparametric problems, it was pointed out in [16] that the worst-case local entropy is of the same order as the global entropy. Therefore a theoretical analysis which relies on global entropy (such as Theorem 4.1) leads to the correct worst-case rate at least in the minimax sense. For parametric problems, at the $O(1/n)$ approximation level, local entropy is constant but the global entropy is $\ln n$. This leads to a $\ln(n)$ difference in the resulting bound.

Although it may not be immediately obvious how to define a localized counterpart of the index of resolvability, we can introduce a correction term which has the same effect. As pointed out earlier, this is essentially the role of the $c_{\rho,n}(\alpha)$ term in Theorem 3.1. We include a simplified version below, which can be obtained by choosing $\alpha = 1/2$ and $\gamma = \rho = 1/\lambda$.

THEOREM 4.2. *Consider the estimator $\hat{k}$ defined in* (6). *Assume that $\lambda > 1$, and let $\rho = 1/\lambda$. Then*

$$\mathbf{E}_X D_\rho^{\mathrm{Re}}(q\|p_{\hat{k}}) \leq \frac{2}{1-\rho} \inf_k \left[ D_{\mathrm{KL}}(q\|p_k) + \frac{\lambda}{n} \ln \frac{\sum_j \pi_j e^{-0.5\rho(1-\rho)n D_\rho^{\mathrm{Re}}(q\|p_j)}}{\pi_k} \right].$$

The bound relies on a localized version of the index of resolvability, with the global entropy $-\ln \pi_k$ replaced by a localized entropy $\ln \sum_j \pi_j \times e^{-0.5\rho(1-\rho)n D_\rho^{\mathrm{Re}}(q\|p_j)} - \ln \pi_k$. Since

$$\ln \sum_j \pi_j e^{-0.5\rho(1-\rho)n D_\rho^{\mathrm{Re}}(q\|p_j)} \leq \ln \sum_j \pi_j = 0,$$

the localized entropy is always smaller than the global entropy. Intuitively, we can see that if $p_j(x)$ is far away from $q(x)$, then $\exp(-\rho(1-\rho)(1-\alpha)n D_\rho^{\mathrm{Re}}(q\|p_j))$ is exponentially small as $n \to \infty$. It follows that the main contribution to the summation in $\sum_j \pi_j e^{-0.5\rho(1-\rho)n D_\rho^{\mathrm{Re}}(q\|p_j)}$ is from terms such that $D_\rho^{\mathrm{Re}}(q\|p_j)$ is small. This is equivalent to a reweighting of the prior $\pi_k$ in such a way that we only count points that are localized within a small $D_\rho^{\mathrm{Re}}$ ball of $q$.

This localization leads to the correct rate of convergence for parametric problems. The effect is similar to using localized entropy in the empirical process analysis. We still consider the same one-dimensional problem discussed at the beginning of the section, with a uniform discretization consisting of $N+1$ points. We will consider the maximum-likelihood estimate. For one-dimensional parametric problems, using the assumption in (8), we have for all $N^2 = O(n)$,

$$\sum_j e^{-\rho(1-\rho)(1-\alpha)n D_\rho^{\mathrm{Re}}(q\|p_j)} \leq \sum_j e^{-\rho(1-\rho)(1-\alpha)n c_1 j^2/N^2} = O(1).$$

Since $\pi_j = 1/(N+1)$, the localized entropy

$$\ln \frac{\sum_j \pi_j e^{-\rho(1-\rho)(1-\alpha)n D_\rho^{\mathrm{Re}}(q\|p_j)}}{\pi_k} = O(1)$$

is a constant when $N = O(n^{1/2})$. Therefore with a discretization size $N = O(n^{1/2})$, Theorem 4.2 implies a convergence rate of the correct order $O(1/n)$.

4.3. *The standard MDL* ($\lambda = 1$). The standard MDL with $\lambda = 1$ in (6) is more complicated to analyze. It is impossible to give a bound similar to Theorem 4.1 that depends only on the index of resolvability. As a matter of fact, no bound was established in [1]. As we will show later, the method can converge very slowly even if the index of resolvability is well behaved.

However, it is possible to obtain bounds in this case under additional assumptions on the rate of decay of the prior $\pi$. The following theorem is a straightforward interpretation of Corollary 3.3, where we consider the family itself as a 0-upper discretization, $\Gamma_i = \{p_i\}$.

THEOREM 4.3. *Consider the estimator defined in* (6) *with* $\lambda = 1$. *For all* $\rho \in (0, 1)$ *and* $\forall \gamma \geq 1$, *we have*

$$\mathbf{E}_X D_\rho^{\mathrm{Re}}(q || p_{\hat{k}}) \leq \frac{\gamma \inf_k [D_{\mathrm{KL}}(q || p_k) + (1/n) \ln(1/\pi_k)]}{\rho(1 - \rho)}$$

$$+ \frac{\gamma - \rho}{\rho(1 - \rho)n} \ln \sum_j \pi_j^{(\gamma - 1)/(\gamma - \rho)}.$$

The above theorem depends only on the index of resolvability and the decay of the prior $\pi$. If $\pi$ has a fast decay in the sense of $\sum_j \pi_j^{(\gamma-1)/(\gamma-\rho)} < +\infty$ and does not change with respect to $n$, then the second term on the right-hand side of Theorem 4.3 is $O(1/n)$. In this case the convergence rate is determined by the index of resolvability. The prior decay condition specified here is rather mild. This implies that the standard MDL is usually Hellinger consistent when used with care.

4.4. *Slow convergence of the standard MDL.* The purpose of this section is to illustrate that the index of resolvability cannot by itself determine the rate of convergence for the standard MDL. We consider a simple example related to the Bayesian inconsistency counterexample given in [2], with an additional randomization argument. Note that due to the randomization, we shall allow two densities in our model class to be identical. It is clear from the construction that this requirement is for convenience only, rather than anything essential.

Given a sample size $n$, consider an integer $m$ such that $m \gg n$. Let the space $\mathcal{X}$ consist of $2m$ points $\{1, \ldots, 2m\}$. Assume that the truth $q$ is the uniform distribution, $q(u) = 1/(2m)$ for $u = 1, \ldots, 2m$.

Consider a density class $\Gamma'$ consisting of all densities $p$ such that either $p(u) = 0$ or $p(u) = 1/m$. That is, a density $p$ in $\Gamma'$ takes the value $1/m$ at $m$ of the $2m$ points, and 0 elsewhere. Now let our model class $\Gamma$ consist of the true density $q$ with prior $1/4$, as well as $2^n$ densities $p_j$ ($j = 1, \ldots, 2^n$) that are randomly and uniformly drawn from $\Gamma'$ (with replacement), where each $p_j$ is given the same prior $3/2^{n+2}$.

We shall show that for a sufficiently large integer $m$, with large probability we will estimate one of the $2^n$ densities from $\Gamma'$ with probability of at least $1 - e^{-1/2}$. Since the index of resolvability is $\ln 4/n$, which is small when $n$ is large, the example implies that the convergence of the standard MDL method cannot be characterized by the index of resolvability alone.

Let $X = \{X_1, \ldots, X_n\}$ be a set of $n$-samples from $q$ and let $\hat{p}$ be the estimator from (6) with $\lambda = 1$ and $\Gamma$ randomly generated above. We would like to estimate $P(\hat{p} = q)$. By construction, $\hat{p} = q$ only when $\prod_{i=1}^{n} p_j(X_i) = 0$ for all $p_j \in \Gamma' \cap \Gamma$. Now pick $m$ large enough such that $(m-n)^n/m^n \geq 0.5$; we have

$$
\begin{aligned}
P(\hat{p} = q) &= P\left( \forall p_j \in \Gamma' \cap \Gamma : \prod_{i=1}^{n} p_j(X_i) = 0 \right) \\
&= \mathbf{E}_X P\left( \forall p_j \in \Gamma' \cap \Gamma : \prod_{i=1}^{n} p_j(X_i) = 0 \Big| X \right) \\
&= \mathbf{E}_X P\left( \prod_{i=1}^{n} p_1(X_i) = 0 \Big| X \right)^{2^n} \\
&= \mathbf{E}_X \left( 1 - \frac{C_{2m-|X|}^{m}}{C_{2m}^{m}} \right)^{2^n} \\
&\leq \mathbf{E}_X \left( 1 - \left( \frac{m-n}{2m} \right)^n \right)^{2^n} \leq (1 - 2^{-(n+1)})^{2^n} \leq e^{-0.5},
\end{aligned}
$$

where $|X|$ denotes the number of distinct elements in $X$. Therefore with a constant probability we have $\hat{p} \neq q$ no matter how large $n$ is.

This example shows that it is impossible to obtain any rate of convergence result using the index of resolvability alone. In order to estimate convergence, it is thus necessary to make additional assumptions, such as the prior decay condition of Theorem 4.3. The randomization used in the construction is not essential. This is because there exists at least one draw (a deterministic configuration) that leads to convergence probability (the probability of correct estimation) at least as large as the expected convergence probability of $e^{-0.5}$ under randomization.

We shall also mention that starting from this example, together with a construction scheme similar to that of the Bayesian inconsistency counterexample in [2], it is not difficult to show that the standard MDL is not Hellinger consistent even when the index of resolvability approaches zero as $n \to \infty$. For simplicity, we skip the detailed construction in this paper.

4.5. *Weak convergence of the standard MDL.* Although Hellinger consistency cannot be obtained for standard MDL based on the index of resolvability alone, it was shown in [1] that as $n \to \infty$, if the index of resolvability approaches zero, then $p_{\hat{k}}$ converges weakly to $q$ in probability (in the sense discussed at the end of Section 3.3). This result is a direct consequence of Theorem 3.2, which we shall restate here.

THEOREM 4.4. *Consider the estimator defined in* (6) *with* $\lambda = 1$. *Then* $\forall f : \mathcal{X} \to [-1, 1]$, *we have*

$$\mathbf{E}_X \left| \mathbf{E}_{p_{\hat{k}}} f(x) - \frac{1}{n} \sum_{i=1}^n f(X_i) \right| \le 2A_n + \sqrt{2A_n},$$

*where* $A_n = \inf_k [D_{\mathrm{KL}}(q||p_k) + \frac{1}{n} \ln \frac{1}{\pi_k}] + \frac{\ln 2}{n}$.

Note that in the sense discussed at the end of Section 3.3, this theorem essentially implies that the standard MDL estimator is weakly consistent (in probability) as long as the index of resolvability approaches zero when $n \to \infty$. Moreover, it establishes a rate of convergence result which depends only on the index of resolvability. This theorem improves the consistency result in [1], where no rate of convergence result was established and $f$ was assumed to be an indicator function.

**5. Bayesian posterior distributions.** Assume we observe $n$-samples $X = \{X_1, \ldots, X_n\} \in \mathcal{X}^n$, independently drawn from the true underlying distribution $Q$ with density $q$. As mentioned earlier, we call any probability density $\hat{w}_X(\theta)$ with respect to $\pi$ that depends on the observation $X$ (and measurable on $\mathcal{X}^n \times \Gamma$) a posterior. For all $\gamma > 0$, we define a generalized Bayesian posterior $\pi_\gamma(\cdot | X)$ with respect to $\pi$ as (also see [15])

$$(9) \qquad \pi_\gamma(\theta | X) = \frac{\prod_{i=1}^n p^\gamma(X_i | \theta)}{\int_\Gamma \prod_{i=1}^n p^\gamma(X_i | \theta) \, d\pi(\theta)}.$$

We call $\pi_\gamma$ the $\gamma$-Bayesian posterior. The standard Bayesian posterior is denoted as $\pi(\cdot | X) = \pi_1(\cdot | X)$.

The key starting point of our analysis is the following simple observation that relates the Bayesian posterior to an instance of information complexity minimization which we have already analyzed in this paper.

PROPOSITION 5.1. *Consider a prior* $\pi$ *and* $\lambda > 0$. *Then*

$$\hat{R}_\lambda(\pi_{1/\lambda}(\cdot | X)) = -\frac{\lambda}{n} \ln \mathbf{E}_\pi \exp\left( \frac{1}{\lambda} \sum_{i=1}^n \ln \frac{p(X_i | \theta)}{q(X_i)} \right) = \inf_w \hat{R}_\lambda(w),$$

*where* $\hat{R}_\lambda(w)$ *is defined in* (2), *and the* inf *on the right-hand side is over all possible densities* $w$ *with respect to the prior* $\pi$.

PROOF.    The first equality follows from simple algebra.

Now let $f(\theta) = \frac{1}{\lambda} \sum_{i=1}^{n} \ln p(X_i|\theta)$ in Proposition 2.1; we obtain

$$-\frac{\lambda}{n} \ln \mathbf{E}_\pi \exp(f(\theta)) \le \inf_w \hat{R}_\lambda(w) \le \hat{R}_\lambda(\pi_{1/\lambda}(\cdot|X)).$$

Combining this with the first equality, we know that equality holds in the above chain of inequalities. This proves the second inequality.    □

The above proposition indicates that the generalized Bayesian posterior can be regarded as a minimum information complexity estimator (1) with $S$ consisting of all possible densities. Therefore results parallel to those of MDL can be obtained.

5.1. *Generalized Bayesian methods.*    Similarly to the index of resolvability complexity measure for MDL, for Bayesian-like methods the corresponding model resolvability, which controls the complexity, becomes the *Bayesian resolvability* defined as

$$\begin{aligned}
(10) \qquad r_{\lambda,n}(q) &= \inf_w \left[ \mathbf{E}_\pi w(\theta) D_{\mathrm{KL}}(q||p(\cdot|\theta)) + \frac{\lambda}{n} D_{\mathrm{KL}}(w\,d\pi||d\pi) \right] \\
&= -\frac{\lambda}{n} \ln \mathbf{E}_\pi e^{-(n/\lambda)D_{\mathrm{KL}}(q||p(\cdot|\theta))}.
\end{aligned}$$

The density that attains the infimum of (10) is given by

$$w(\theta) \propto \exp\left[ -\frac{n}{\lambda} D_{\mathrm{KL}}(q||p(\cdot||\theta)) \right].$$

The following proposition gives a simple and intuitive estimate of the Bayesian index of resolvability. This bound implies that the Bayesian resolvability can be estimated using local properties of the prior $\pi$ around the true density $q$. The quantity is small as long as there is a positive prior mass in a small KL-ball around the truth $q$.

PROPOSITION 5.2.    *The Bayesian resolvability defined in (10) can be bounded as*

$$r_{\lambda,n}(q) \le \inf_{\varepsilon > 0} \left[ \varepsilon - \frac{\lambda}{n} \ln \pi(\{p \in \Gamma : D_{\mathrm{KL}}(q||p) \le \varepsilon\}) \right].$$

PROOF.    For all $\varepsilon > 0$, we simply note that $\mathbf{E}_\pi e^{-(n/\lambda)D_{\mathrm{KL}}(q||p(\cdot|\theta))} \ge e^{-(n/\lambda)\varepsilon} \times \pi(\{p \in \Gamma : D_{\mathrm{KL}}(q||p) \le \varepsilon\})$. Now taking the logarithm and using (10), we obtain the desired inequality.    □

The following bound is a direct consequence of Corollary 3.2.

Theorem 5.1. *Consider the generalized Bayesian posterior $\pi_{1/\lambda}(\theta|X)$ defined in (9) with $\lambda > 1$. Then $\forall \rho \in (0, 1/\lambda]$*

$$\mathbf{E}_X \mathbf{E}_\pi \pi_{1/\lambda}(\theta|X) D_\rho^{\mathrm{Re}}(q||p(\cdot|\theta)) \leq -\frac{\lambda}{\rho(\lambda-1)n} \ln \mathbf{E}_\pi \exp\left(-\frac{n}{\lambda} D_{\mathrm{KL}}(q||p(\cdot|\theta))\right).$$

The above theorem gives a general convergence bound on the $\gamma$-Bayesian method with $\gamma < 1$, depending only on the globally defined Bayesian resolvability. Note that similarly to Theorem 4.2 for the MDL case, a bound using a localized Bayesian resolvability can also be obtained.

Theorem 5.1 immediately implies the concentration of a generalized Bayesian posterior. Define the posterior mass outside an $\varepsilon$ $D_\rho^{\mathrm{Re}}$-ball around $q$ as

$$\pi_{1/\lambda}(\{p \in \Gamma : D_\rho^{\mathrm{Re}}(q||p) \geq \varepsilon\}|X).$$

Using the bound in Theorem 5.1 and Proposition 5.2, we can show that with large probability, the generalized Bayesian posterior outside a $D_\rho^{\mathrm{Re}}$-ball of size $O(\varepsilon)$ is exponentially small when $\varepsilon \gg \varepsilon_{\pi,n}$. However, the average performance bound in Theorem 5.1 is not refined enough to yield exponential tail probability directly under the prior $\pi$. In order to obtain the correct behavior, we shall thus consider a prior $\pi'$ related to $\pi$ which is more heavily concentrated on distributions that are far away from $q$. We choose $\pi'$ for which Theorem 5.1 can be used to obtain a constant probability of posterior concentration. We then translate the concentration of posterior with respect to $\pi'$ to a concentration result with respect to $\pi$.

Corollary 5.1. *Let $\lambda > 1$ and $\rho \in (0, 1/\lambda]$. Then for all $t \geq 0$ and $\delta \in (0, 1)$, with probability at least $1 - \delta$,*

$$\pi_{1/\lambda}\left(\left\{p \in \Gamma : D_\rho^{\mathrm{Re}}(q||p) \geq \frac{4\varepsilon_{\pi,n} + 2t}{\rho(\lambda-1)\delta}\right\}\Big|X\right) \leq \frac{1}{1 + e^{nt/\lambda}},$$

*where the critical prior-mass radius $\varepsilon_{\pi,n} = \inf\{\varepsilon : \varepsilon \geq -\frac{\lambda}{n} \ln \pi(\{p \in \Gamma : D_{\mathrm{KL}}(q||p) \leq \varepsilon\})\}$.*

Proof. Let $\varepsilon_t = 2(2\varepsilon_{\pi,n} + t)/(\rho(\lambda-1)\delta)$, $\Gamma_1 = \{p \in \Gamma : D_\rho^{\mathrm{Re}}(q||p) < \varepsilon_t\}$ and $\Gamma_2 = \{p \in \Gamma : D_\rho^{\mathrm{Re}}(q||p) \geq \varepsilon_t\}$. We let $a = e^{-nt/\lambda}$ and define $\pi'(\theta) = a\pi(\theta)C$ when $\theta \in \Gamma_1$ and $\pi'(\theta) = \pi(\theta)C$ when $\theta \in \Gamma_2$, where the normalization constant $C = (a\pi(\Gamma_1) + \pi(\Gamma_2))^{-1} \in [1, 1/a]$.

Now apply Theorem 5.1 and Proposition 5.2 with the prior $\pi'$. We obtain (using the Markov inequality)

$$\mathbf{E}_X \pi'_{1/\lambda}(\Gamma_2|X) \varepsilon_t \leq \mathbf{E}_X \mathbf{E}_{\pi'} \pi'_{1/\lambda}(\theta|X) D_\rho^{\mathrm{Re}}(q||p)$$

$$\leq -\frac{\lambda}{\rho(\lambda-1)n} \ln \mathbf{E}_{\pi'} \exp\left(-\frac{n}{\lambda} D_{\mathrm{KL}}(q||p(\cdot|\theta))\right)$$

$$\leq -\frac{\lambda}{\rho(\lambda-1)n}\left[\ln a + \ln \mathbf{E}_\pi \exp\left(-\frac{n}{\lambda}D_{\mathrm{KL}}(q||p(\cdot|\theta))\right)\right]$$

$$\leq \frac{1}{\rho(\lambda-1)}\left[t + \varepsilon_{\pi,n} - \frac{\lambda}{n}\ln \pi(\{p \in \Gamma : D_{\mathrm{KL}}(q||p) \leq \varepsilon_{\pi,n}\})\right]$$

$$\leq \frac{2\varepsilon_{\pi,n}+t}{\rho(\lambda-1)}.$$

In the above derivation, the first inequality is the Markov inequality; the second inequality is from Theorem 5.1; the third inequality follows from $\pi' \geq a\pi$; the fourth inequality follows from Proposition 5.2; the final inequality uses the definition of $\varepsilon_{\pi,n}$.

Now we can divide both sides by $\varepsilon_t$, and obtain with probability $1-\delta$ that $\pi'_{1/\lambda}(\Gamma_2|X) \leq 0.5$. By construction, $\pi'_{1/\lambda}(\Gamma_2|X) = \pi_{1/\lambda}(\Gamma_2|X)/(a(1-\pi_{1/\lambda}(\Gamma_2|X)) + \pi_{1/\lambda}(\Gamma_2|X))$. We can solve for $\pi_{1/\lambda}(\Gamma_2|X)$ as $\pi_{1/\lambda}(\Gamma_2|X) = a\pi'_{1/\lambda}(\Gamma_2|X)/(1-(1-a)\pi'_{1/\lambda}(\Gamma_2|X)) \leq a/(1+a)$. □

From the bound, we can see that with large probability the posterior probability outside a $D_\rho^{\mathrm{Re}}$-ball with large distance $t$ decays exponentially in $nt$ and is independent of the complexity of the prior (as long as $t$ is larger than the scale of the critical radius $\varepsilon_{\pi,n}$). As we will see later, the same is true for the standard Bayesian posterior distributions.

5.2. *The standard Bayesian method.* For the standard Bayesian posterior distribution, it is impossible to bound its convergence using only the Bayesian resolvability. The reason is the same as in the MDL case. In fact, it is immediately obvious that the example for MDL can also be applied here. Also see [2] for a related example.

Therefore in order to obtain a rate of convergence (and concentration) for the standard Bayesian method, additional assumptions are necessary. Similarly to Theorem 4.3, bounds using upper-bracketing radius can be easily obtained from Corollary 3.3.

THEOREM 5.2. *Consider the Bayesian posterior $\pi(\cdot|X) = \pi_1(\cdot|X)$ defined in (9). Consider an arbitrary cover $\{\Gamma_j\}$ of $\Gamma$. Then $\forall \rho \in (0,1)$ and $\gamma \geq 1$, we have*

$$\mathbf{E}_X \mathbf{E}_\pi \pi(\theta|X) D_\rho^{\mathrm{Re}}(q||p(\cdot|\theta))$$

$$\leq \frac{\gamma \ln \mathbf{E}_\pi e^{-nD_{\mathrm{KL}}(q||p(\cdot|\theta))}}{\rho(\rho-1)n}$$

$$+ \frac{\gamma-\rho}{\rho(1-\rho)n}\ln \sum_j \pi(\Gamma_j)^{(\gamma-1)/(\gamma-\rho)}(1+r_{\mathrm{ub}}(\Gamma_j))^n.$$

For all $\varepsilon > 0$, consider an $\varepsilon$-upper discretization $\{\Gamma_j^\varepsilon\}$ of $\Gamma$. We obtain from Theorem 5.2,

$$\mathbf{E}_X \mathbf{E}_\pi \pi(\theta|X) D_\rho^{\mathrm{Re}}(q||p(\cdot|\theta))$$

$$\leq \frac{\gamma \ln \mathbf{E}_\pi e^{-n D_{\mathrm{KL}}(q||p(\cdot|\theta))}}{\rho(\rho-1)n}$$

$$+ \frac{\gamma - \rho}{\rho(1-\rho)} \inf_{\varepsilon > 0} \left[ \frac{\ln \sum_j \pi(\Gamma_j^\varepsilon)^{(\gamma-1)/(\gamma-\rho)}}{n} + \ln(1+\varepsilon) \right].$$

In particular, let $\gamma \to 1$. We have

$$\mathbf{E}_X \mathbf{E}_\pi \pi(\theta|X) D_\rho^{\mathrm{Re}}(q||p(\cdot|\theta))$$

$$\leq \frac{\ln \mathbf{E}_\pi e^{-n D_{\mathrm{KL}}(q||p(\cdot|\theta))}}{\rho(\rho-1)n} + \frac{1}{\rho} \inf_{\varepsilon > 0} \left[ \frac{\ln N_{\mathrm{ub}}(\Gamma, \varepsilon)}{n} + \ln(1+\varepsilon) \right],$$

where $N_{\mathrm{ub}}(\Gamma, \varepsilon)$ is the $\varepsilon$-upper-bracketing covering number of $\Gamma$.

Similarly to Corollary 5.1, we obtain the following concentration result for the standard Bayesian posterior distribution from Theorem 5.2.

COROLLARY 5.2. Let $\varepsilon_{\pi, n} = \inf\{\varepsilon : \varepsilon \geq -\frac{1}{n} \ln \pi(\{p \in \Gamma : D_{\mathrm{KL}}(q||p) \leq \varepsilon\})\}$ be the critical prior-mass radius. Let $\rho \in (0, 1)$. For all $s \in [0, 1]$, let

$$\varepsilon_{\mathrm{upper}, n}(s) = \frac{1}{n} \inf_{\{\Gamma_j\}} \ln \sum_j \pi(\Gamma_j)^s (1 + r_{\mathrm{ub}}(\Gamma_j))^n$$

be the critical upper-bracketing radius with coefficient $s$, where $\{\Gamma_j\}$ denotes an arbitrary covering of $\Gamma$. Now $\forall \rho \in (0, 1)$ and $\gamma \geq 1$, let

$$\varepsilon_n = 2\gamma \varepsilon_{\pi, n} + (\gamma - \rho) \varepsilon_{\mathrm{upper}, n}((\gamma - 1)/(\gamma - \rho)).$$

We have for all $t \geq 0$ and $\delta \in (0, 1)$, with probability at least $1 - \delta$,

$$\pi\left( \left\{ p \in \Gamma : D_\rho^{\mathrm{Re}}(q||p) \geq \frac{2\varepsilon_n + (4\gamma - 2)t}{\rho(1-\rho)\delta} \right\} \Big| X \right) \leq \frac{1}{1 + e^{nt}}.$$

PROOF. The proof is similar to that of Corollary 5.1. We let $\varepsilon_t = (2\varepsilon_n + (4\gamma - 2)t)/((\rho - \rho^2)\delta)$. Define $\Gamma_1 = \{p \in \Gamma : D_\rho^{\mathrm{Re}}(q||p) < \varepsilon_t\}$ and $\Gamma_2 = \{p \in \Gamma : D_\rho^{\mathrm{Re}}(q||p) \geq \varepsilon_t\}$. We let $a = e^{-nt}$ and define $\pi'(\theta) = a\pi(\theta)C$ when $\theta \in \Gamma_1$ and $\pi'(\theta) = \pi(\theta)C$ when $\theta \in \Gamma_2$, where the normalization constant $C = (a\pi(\Gamma_1) + \pi(\Gamma_2))^{-1} \in [1, 1/a]$.

Using Proposition 5.2 and the assumption of the theorem, we obtain

$$\frac{\gamma \ln \mathbf{E}_{\pi'} e^{-n D_{\mathrm{KL}}(q||p(\cdot|\theta))}}{\rho(\rho-1)n}$$

$$+ \frac{\gamma - \rho}{\rho(1-\rho)n} \inf_{\{\Gamma_j\}} \ln \sum_j \pi'(\Gamma_j)^{(\gamma-1)/(\gamma-\rho)} (1 + r_{\mathrm{ub}}(\Gamma_j))^n$$

$$\leq \frac{\gamma t + (\gamma/n) \ln \mathbf{E}_\pi e^{-n D_{\mathrm{KL}}(q\|p(\cdot|\theta))}}{\rho(1-\rho)}$$

$$+ \frac{\gamma - \rho}{\rho(1-\rho)} \left[ \frac{(\gamma-1)t}{\gamma - \rho} + \varepsilon_{\mathrm{upper},n} \left( \frac{\gamma-1}{\gamma-\rho} \right) \right]$$

$$\leq \frac{(2\gamma-1)t + \varepsilon_n}{\rho(1-\rho)}.$$

In the first inequality, we have used the fact that $a\pi(\theta) \leq \pi'(\theta) \leq \pi(\theta)/a$. Similarly to the proof of Corollary 5.1, we can use Markov inequality to obtain $\pi'(\Gamma_2|X) \leq 0.5$ with probability $1 - \delta$. This leads to the desired bound for $\pi(\Gamma_2|X) = a\pi'_{1/\lambda}(\Gamma_2|X)/(1 - (1-a)\pi'_{1/\lambda}(\Gamma_2|X))$.  □

In this theorem, we can use the estimate

$$\varepsilon_{\mathrm{upper},n}(s) \leq \inf_{\varepsilon>0} \left[ \frac{1}{n} \ln N_{\mathrm{ub}}(\Gamma, \varepsilon) + \ln(1+\varepsilon) \right],$$

where $N_{\mathrm{ub}}(\Gamma, \varepsilon)$ is the upper-bracketing covering number of $\Gamma$ at scale $\varepsilon$. The result implies that if the critical upper-bracketing radius $\varepsilon_{\mathrm{upper},n}$ is at the same (or smaller) order of the critical prior-mass radius $\varepsilon_{\pi,n}$, then with large probability, the standard Bayesian posterior distribution will concentrate in a $D_\rho^{\mathrm{Re}}$-ball of size $\varepsilon_{\pi,n}$. In this case, the standard Bayesian posterior has the same rate of convergence when compared with the generalized Bayesian posterior with $\lambda > 1$. However, if $\varepsilon_{\mathrm{upper},n}$ is large, then the standard Bayesian method may fail to concentrate in a small $D_\rho^{\mathrm{Re}}$-ball around the truth $q$, even when the critical prior radius $\varepsilon_{\pi,n}$ is small. This can be easily seen from the same counterexample used to illustrate the slow convergence of the standard MDL.

Although the standard Bayesian posterior distribution may not concentrate even when $\varepsilon_{\pi,n}$ is small, Theorem 3.2 implies that the Bayesian density estimator $\mathbf{E}_\pi \pi(\theta|X)p(\cdot|X)$ is close to $q$ in the sense of weak convergence.

The consistency theorem given in [2] also relies on the upper covering number $N_{\mathrm{ub}}(\Gamma, \varepsilon)$. However, no convergence rate was established. Therefore Corollary 5.2 in some sense can be regarded as a refinement of their analysis using their covering definition. Other kinds of covering numbers (e.g., Hellinger covering) can also be used in convergence analysis of nonparametric Bayesian methods. For example, some different definitions can be found in [4] and [12].

The convergence analysis in [12] employed techniques from empirical processes, which can possibly lead to suboptimal convergence rates when the covering number grows relatively fast as the scale $\varepsilon \to 0$. We shall focus on [4], which employed techniques from hypothesis testing in [6]. The resulting convergence theorem from their analysis cannot be as simply stated as those in this paper. Moreover, some of their conditions can be relaxed. Using techniques of this paper, we can obtain the following result. The proof, which requires two additional lemmas, is left to the Appendix.

**THEOREM 5.3.** *Consider a partition of $\Gamma$ as the union of countably many disjoint measurable sets $\Gamma_j$ $(j = 1, \ldots)$. Then $\forall \rho \in (0, 1)$ and $\gamma \geq 1$*

$$\mathbf{E}_X \sum_j \pi(\Gamma_j | X) \inf_{p \in \mathrm{co}(\Gamma_j)} D_\rho^{\mathrm{Re}}(q \| p)$$

$$\leq \frac{(\gamma - \rho) \ln \sum_j \pi(\Gamma_j)^{(\gamma - 1)/(\gamma - \rho)} - \gamma \ln \sum_j \pi(\Gamma_j) e^{-n \sup_{p \in \mathrm{co}(\Gamma_j)} D_{\mathrm{KL}}(q \| p)}}{\rho(1 - \rho)n},$$

*where $\mathrm{co}(\Gamma_j)$ is the convex hull of densities in $\Gamma_j$, $\pi(\Gamma_j) = \int_{\Gamma_j} d\pi(\theta)$ is the prior probability of $\Gamma_j$ and $\pi(\Gamma_j | X) = \int_{\Gamma_j} \prod_{i=1}^n p(X_i | \theta) \, d\pi(\theta) / \int_\Gamma \prod_{i=1}^n p(X_i | \theta) \, d\pi(\theta)$ is the Bayesian posterior probability of $\Gamma_j$.*

An immediate consequence of the above theorem is a result on the concentration of Bayesian posterior distributions that refines some aspects of the main result in [4]. It also complements the upper-bracketing radius-based bound in Corollary 5.2. For simplicity, we only state a version for $\rho$-divergence so that the result is directly comparable to that of [4]. A similar bound can be stated for Rényi entropy.

**COROLLARY 5.3.** *Let $\varepsilon_{\pi,n} = \inf\{\varepsilon : \varepsilon \geq -\frac{1}{n} \ln \pi(\{p \in \Gamma : D_{\mathrm{KL}}(q \| p) \leq \varepsilon\})\}$. Given $\rho \in (0, 1)$, we assume that $\forall \varepsilon > 0$, $\{p \in \Gamma : D_\rho(q \| p) \geq \varepsilon\}$ can be covered by the union of measurable sets $\Gamma_j^\varepsilon$ $(j = 1, \ldots)$ such that $\inf\{D_\rho(q \| p) : p \in \bigcup_j \mathrm{co}(\Gamma_j^\varepsilon)\} \geq \varepsilon/2$. For all $s \in [0, 1]$, let*

$$\varepsilon_{\mathrm{conv},n}(s) = \sup\left\{\varepsilon_0 : \varepsilon_0 < \frac{1}{n} \sup_{\varepsilon \geq \varepsilon_0} \inf_{\{\Gamma_j^\varepsilon\}} \ln\left(\sum_j \pi(\Gamma_j^\varepsilon)^s + 2\right)\right\}$$

*be the critical convex-cover radius. Now $\forall \gamma \geq 1$ let*

$$\varepsilon_n = 2\gamma \varepsilon_{\pi,n} + (\gamma - \rho) \varepsilon_{\mathrm{conv},n}((\gamma - 1)/(\gamma - \rho)).$$

*For all $t \geq 0$ and $\delta \in (0, 1)$, with probability at least $1 - \delta$,*

$$\pi\left(\left\{p \in \Gamma : D_\rho(q \| p) \geq \frac{4\varepsilon_n + (8\gamma - 4)t}{\rho(1 - \rho)\delta}\right\} \Big| X\right) \leq \frac{1}{1 + e^{nt}}.$$

**PROOF.** Let $\varepsilon_t = 4(\varepsilon_n + (2\gamma - 1)t)/(\rho(1 - \rho)\delta)$. Similarly to the proof of Corollary 5.1, we define $\Gamma_1 = \{p \in \Gamma : D_\rho(q \| p) < \varepsilon_t\}$, $\Gamma_2 = \Gamma - \Gamma_1$. We let $a = e^{-nt}$ and define $\pi'(\theta) = a\pi(\theta)C$ when $\theta \in \Gamma_1$ and $\pi'(\theta) = \pi(\theta)C$ when $\theta \in \Gamma_2$, where the normalization constant $C = (a\pi(\Gamma_1) + \pi(\Gamma_2))^{-1} \in [1, 1/a]$.

Let $\Gamma'_0 = \{p \in \Gamma : D_{\mathrm{KL}}(q||p) < \varepsilon_{\pi,n}\}$. Since $D_{\mathrm{KL}}(q||p) = D_0(q||p)$ and $\varepsilon_t \geq \varepsilon_{\pi,n}/\min(\rho, 1-\rho)$, we know from Proposition 3.1 that $\Gamma'_0 \subset \Gamma_1$. Let $\Gamma'_{-1} = \Gamma_1 - \Gamma'_0$. By assumption, it is clear that $\Gamma_2$ can be partitioned into the union of disjoint measurable sets $\{\Gamma'_j\}$ ($j \geq 1$) such that $\Gamma'_j \subset \Gamma_j^{\varepsilon_t}$ and $\inf_{p \in \bigcup_{j \geq 1} \mathrm{co}(\Gamma'_j)} D_\rho(q||p) \geq \varepsilon_t/2$. For this partition, we have

$$\mathbf{E}_X \pi'(\Gamma_2|X)\varepsilon_t/2 \leq \mathbf{E}_X \sum_{j \geq -1} \pi'(\Gamma'_j|X) \inf_{p \in \mathrm{co}(\Gamma'_j)} D_\rho(q||p).$$

Note that

$$\ln \sum_{j \geq -1} \pi'(\Gamma'_j)^{(\gamma-1)/(\gamma-\rho)} \leq -\frac{\gamma-1}{\gamma-\rho} \ln a + \ln\left[\sum_{j \geq 1} \pi(\Gamma_j^{\varepsilon_t})^{(\gamma-1)/(\gamma-\rho)} + 2\right]$$

$$\leq -\frac{\gamma-1}{\gamma-\rho} \ln a + n\varepsilon_{\mathrm{conv},n}$$

and

$$-\ln \sum_{j \geq -1} \pi'(\Gamma'_j) e^{-n\sup_{p \in \mathrm{co}(\Gamma'_j)} D_{\mathrm{KL}}(q||p)} \leq n \sup_{p \in \mathrm{co}(\Gamma'_0)} D_{\mathrm{KL}}(q||p) - \ln \pi'(\Gamma'_0)$$

$$\leq 2n\varepsilon_{\pi,n} + nt.$$

Combining the above estimates, and plugging them into Theorem 5.3, we obtain

$$\mathbf{E}_X \pi'(\Gamma_2|X) \leq \frac{(\gamma-\rho)(-(\ln a)(\gamma-1)/(\gamma-\rho) + n\varepsilon_{\mathrm{conv},n}) + \gamma(2n\varepsilon_{\pi,n} + nt)}{\rho(1-\rho)n\varepsilon_t/2}$$

$$= 0.5\delta.$$

Therefore $\pi'(\Gamma_2|X) \leq 0.5$ with probability $1 - \delta$. The desired bound for $\pi(\Gamma_2|X)$ can be obtained from $\pi(\Gamma_2|X) = a\pi'_{1/\lambda}(\Gamma_2|X)/(1 - (1-a)\pi'_{1/\lambda}(\Gamma_2|X))$. $\square$

If we can cover $\{p \in \Gamma : D_\rho(q||p) \geq \varepsilon\}$ by $N_\varepsilon$ convex measurable sets $\Gamma_j^\varepsilon$ ($j = 1, \ldots, N_\varepsilon$) such that $\inf\{D_\rho(q||p) : p \in \bigcup_j \Gamma_j^\varepsilon\} \geq \varepsilon/2$, then we may take $\gamma = 1$ in Corollary 5.3 with $\varepsilon_{\mathrm{conv},n}$ defined as

$$\varepsilon_{\mathrm{conv},n} = \sup\left\{\varepsilon_0 : \varepsilon_0 < \frac{1}{n} \ln\left(\sup_{\varepsilon \geq \varepsilon_0} N_\varepsilon + 2\right)\right\}.$$

Clearly if $\frac{1}{n} \ln N_\varepsilon = O(\varepsilon_{\pi,n})$ for some $\varepsilon = O(\varepsilon_{\pi,n})$, then with large probability Bayesian posterior distributions concentrate on a $D_\rho$-ball of size $O(\varepsilon_{\pi,n})$ around $q$. Note that this result relaxes a condition of [4], where our definition of $\varepsilon_{\pi,n}$ was replaced by possibly smaller balls $\{p \in \Gamma : D_{\mathrm{KL}}(q||p) \leq$ $\varepsilon, \mathbf{E}_q \ln(\frac{q}{p})^2 \le \varepsilon\}$. Moreover, their covering definition $N_\varepsilon$ does not apply to arbitrary convex covering sets directly (although it is not difficult to modify their proof to deal with this case), and their result does not directly handle noncompact families where $N_\varepsilon = \infty$ (which can be directly handled by our result with $\gamma > 1$).

It is worth mentioning that for practical purposes, the balls $\{p \in \Gamma : D_{\mathrm{KL}}(q\|p) \le \varepsilon, \mathbf{E}_q \ln(\frac{q}{p})^2 \le \varepsilon\}$ and $\{p \in \Gamma : D_{\mathrm{KL}}(q\|p) \le \varepsilon\}$ are usually of comparable size. Therefore relaxing this condition may not always lead to significant practical advantages. However, it is possible to construct examples such that this refinement makes a difference. For example, consider the discrete family $\Gamma = \{p_j\}$ $(j \ge 1)$ with prior $\pi_j = 1/j(j+1)$. Assume that the truth $q(x)$ is the uniform distribution on $[0, 1]$, and $p_j(x) = 2^{-j}$ when $x \in [0, j^{-2}/2]$ and $p_j(x) = (j^2 - 2^{-j-1})/(j^2 - 0.5)$ otherwise. It is clear that $\mathbf{E}_q \ln(\frac{q}{p_j})^2 \ge 0.5 \ln 4$, while $\lim_{j\to\infty} D_{\mathrm{KL}}(q\|p_j) = 0$. Therefore the result in [4] cannot be applied, while Corollary 5.3 implies that the posterior distribution is consistent in this example.

Applications of convergence results similar to Corollary 5.2 and Corollary 5.3 can be found in [4] and [12]. It is also useful to note that Corollary 5.1 requires less assumptions to achieve good convergence rates, implying that generalized Bayesian methods are more stable than the standard Bayesian method. This fact has also been observed in [15].

## 6. Discussion.
This paper studies certain randomized (and deterministic) density estimation methods which we call information complexity minimization. We introduced a general KL-entropy based convergence analysis, and demonstrated that this approach can lead to simplified and improved convergence results for MDL and Bayesian posterior distributions.

An important observation from our study is that generalized information complexity minimization methods with regularization parameter $\lambda > 1$ are more robust than the corresponding standard methods with $\lambda = 1$. That is, their convergence behavior is completely determined by the local prior density around the true distribution measured by the model resolvability $\inf_{w\in S} R_\lambda(w)$. For MDL, this quantity (index of resolvability) is well behaved if we put a not too small prior mass at a density that is close to the truth $q$. For the Bayesian posterior, this quantity (Bayesian resolvability) is well behaved if we put a not too small prior mass in a small KL-ball around $q$. We have also demonstrated through an example that the standard MDL (and Bayesian posterior) does not have this desirable property. That is, even if we can guess the true density by putting a relatively large prior mass at the true density $q$, we may not be able to estimate $q$ very well as long as there exists a bad (random) prior structure even at places very far from the truth $q$.

Therefore, although the standard Bayesian method is "optimal" in a certain averaging sense, its behavior is heavily dependent on the regularity of the prior distribution globally. Intuitively, the standard Bayesian method can put too much emphasis on the difficult part of the prior distribution, which degrades the estimation quality in the easier part in which we are actually more interested. Therefore even if one is able to guess the true distribution by putting a large prior mass around its neighborhood, the Bayesian method can still behave poorly if one accidentally makes bad choices elsewhere. This implies that unless one completely understands the impact of the prior, it is much safer to use a generalized Bayesian method with $\lambda > 1$.

## APPENDIX

**A.1. Lower bounds of $\mathbf{E}_X \hat{R}_{\lambda'}(\hat{w}_X^S)$.** In order to apply Theorem 3.1, we shall bound the quantity $\mathbf{E}_X \hat{R}_{\lambda'}(\hat{w}_X^S)$ from below.

LEMMA A.1.    *For all $\lambda' \geq 1$, $\mathbf{E}_X \hat{R}_{\lambda'}(\hat{w}_X^S) \geq -\frac{\lambda'}{n} \ln \mathbf{E}_\pi \mathbf{E}_q^n (\frac{p(x|\theta)}{q(x)})^{1/\lambda'} \geq 0$.*

PROOF.    The convex duality in Proposition 2.1 with $f(x) = -\frac{1}{\lambda'} \sum_{i=1}^n \ln \frac{q(X_i)}{p(X_i|\theta)}$ implies

$$\hat{R}_{\lambda'}(\hat{w}_X^S) \geq -\frac{\lambda'}{n} \ln \mathbf{E}_\pi \exp\left(-\frac{1}{\lambda'} \sum_{i=1}^n \ln \frac{q(X_i)}{p(X_i|\theta)}\right).$$

Now by taking expectation and using Jensen's inequality with the convex function $\psi(x) = -\ln(x)$, we obtain

$$\mathbf{E}_X \hat{R}_{\lambda'}(\hat{w}_X^S) \geq -\frac{\lambda'}{n} \ln \mathbf{E}_X \mathbf{E}_\pi \exp\left(-\frac{1}{\lambda'} \sum_{i=1}^n \ln \frac{q(X_i)}{p(X_i|\theta)}\right)$$

$$= -\frac{\lambda'}{n} \ln \mathbf{E}_\pi \mathbf{E}_q^n \left(\frac{p(x|\theta)}{q(x)}\right)^{1/\lambda'} \geq 0,$$

which proves the lemma.    □

LEMMA A.2.    *Consider an arbitrary cover $\{\Gamma_j\}$ of $\Gamma$. The following inequality is valid $\forall \lambda' \in [0, 1]$:*

$$\mathbf{E}_X \hat{R}_{\lambda'}(\hat{w}_X^S) \geq -\frac{1}{n} \ln \sum_j \pi(\Gamma_j)^{\lambda'} (1 + r_{\mathrm{ub}}(\Gamma_j))^n,$$

*where $r_{\mathrm{ub}}$ is the upper-bracketing radius in Definition 3.2.*

PROOF. The proof is similar to that of Lemma A.1, but with a slightly different estimate. We again start with the inequality

$$\hat{R}_{\lambda'}(\hat{w}_X^S) \geq -\frac{\lambda'}{n}\ln \mathbf{E}_\pi \exp\left(-\frac{1}{\lambda'}\sum_{i=1}^n \ln\frac{q(X_i)}{p(X_i|\theta)}\right).$$

Taking expectation and using Jensen's inequality with the convex function $\psi(x) = -\ln(x)$, we obtain

$$-\mathbf{E}_X \hat{R}_{\lambda'}(\hat{w}_X^S) \leq \frac{1}{n}\ln \mathbf{E}_X \mathbf{E}_\pi^{\lambda'} \exp\left(-\frac{1}{\lambda'}\sum_{i=1}^n \ln\frac{q(X_i)}{p(X_i|\theta)}\right)$$

$$\leq \frac{1}{n}\ln \mathbf{E}_X\left[\sum_j \pi(\Gamma_j)\exp\left(-\frac{1}{\lambda'}\sum_{i=1}^n \ln\frac{q(X_i)}{\sup_{\theta\in\Gamma_j}p(X_i|\theta)}\right)\right]^{\lambda'}$$

$$\leq \frac{1}{n}\ln \mathbf{E}_X\left[\sum_j \pi(\Gamma_j)^{\lambda'}\exp\left(-\sum_{i=1}^n \ln\frac{q(X_i)}{\sup_{\theta\in\Gamma_j}p(X_i|\theta)}\right)\right]$$

$$= \frac{1}{n}\ln\left[\sum_j \pi(\Gamma_j)^{\lambda'}\mathbf{E}_X \prod_{i=1}^n \frac{\sup_{\theta\in\Gamma_j}p(X_i|\theta)}{q(X_i)}\right]$$

$$= \frac{1}{n}\ln\left[\sum_j \pi(\Gamma_j)^{\lambda'}(1+r_{\mathrm{ub}}(\Gamma_j))^n\right].$$

The third inequality follows from the fact that $\forall \lambda' \in [0,1]$ and positive numbers $\{a_j\}$, $(\sum_j a_j)^{\lambda'} \leq \sum_j a_j^{\lambda'}$. $\quad\square$

## A.2. Proof of Theorem 5.3. The proof requires two lemmas.

LEMMA A.3. *Consider a partition of $\Gamma$ as the union of countably many disjoint measurable sets $\Gamma_j$ $(j = 1,\ldots)$. Let*

$$q(X) = \prod_{i=1}^n q(X_i), \qquad p_j(X) = \frac{1}{\pi(\Gamma_j)}\int_{\Gamma_j}\prod_{i=1}^n p(X_i|\theta)\,d\pi(\theta).$$

*Then we have $\forall \rho \in (0,1)$ and $\gamma \geq 1$,*

$$\mathbf{E}_X \sum_j \pi(\Gamma_j|X)D_\rho^{\mathrm{Re}}(q(X')||p_j(X'))$$

$$\leq \frac{(\gamma-\rho)\ln\sum_j \pi(\Gamma_j)^{(\gamma-1)/(\gamma-\rho)} - \gamma\ln\sum_j \pi(\Gamma_j)e^{-D_{\mathrm{KL}}(q(X')||p_j(X'))}}{\rho(1-\rho)},$$

*where $X', X \in \mathcal{X}^n$, $q(X) = \prod_{i=1}^n q(X_i)$ is the true density of $X$ and*

$$p_j(X) = \frac{1}{\pi(\Gamma_j)}\int_{\Gamma_j}\prod_{i=1}^n p(X_i|\theta)\,d\pi(\theta)$$

*is the mixture density over* $\Gamma_j$ *under* $\pi$.

PROOF. We shall apply Corollary 3.3 with a slightly different interpretation. Instead of considering $X$ as $n$ independent samples $X_i$ as before, we simply regard it as one random variable by itself. Consider the family $\Gamma'$ which consists of discrete densities $p_j(X)$, with prior $\pi_j = \pi(\Gamma_j)$. This discretization itself can be regarded as a 0-upper discretization of $\Gamma'$. Also, given $X$, it is easy to see that the Bayesian posterior on $\Gamma'$ with respect to $\{\pi_j\}$ is $\hat{\pi}_j = \pi(\Gamma_j | X)$. We can thus apply Corollary 3.3 on $\Gamma'$, which leads to the stated bound [with the help of (10)]. □

In order to apply the above lemma, we also need to simplify $D_\rho^{\mathrm{Re}}(q(X') \| p_j(X'))$ and $D_{\mathrm{KL}}(q(X') \| p_j(X'))$.

LEMMA A.4. *We have the bounds*

$$\inf_{p \in \mathrm{co}(\Gamma_j)} D_\rho^{\mathrm{Re}}(q(X_1) \| p(X_1)) \leq \frac{D_\rho^{\mathrm{Re}}(q(X) \| p_j(X))}{n}$$

$$\leq \sup_{p \in \mathrm{co}(\Gamma_j)} D_\rho^{\mathrm{Re}}(q(X_1) \| p(X_1))$$

*and*

$$\inf_{p \in \mathrm{co}(\Gamma_j)} D_{\mathrm{KL}}(q(X_1) \| p(X_1)) \leq \frac{D_{\mathrm{KL}}(q(X) \| p_j(X))}{n}$$

$$\leq \sup_{p \in \mathrm{co}(\Gamma_j)} D_{\mathrm{KL}}(q(X_1) \| p(X_1)).$$

PROOF. Since $D_{\mathrm{KL}}(q \| p) = \lim_{\rho \to 0^+} D_\rho^{\mathrm{Re}}(q \| p)$, we only need to prove the first two inequalities. The proof is essentially the same as that of Lemma 4 on page 478 of [6], which dealt with the existence of tests under the Hellinger distance. We include it here for completeness.

We shall only prove the first half of the first two inequalities (the second half has an identical proof) and we shall prove the claim by induction. If $n = 1$, then since $p_j(X) \in \mathrm{co}(\Gamma_j)$ the claim holds trivially. Now assume that the claim holds for $n = k$. For $n = k + 1$, if we let

$$w(\theta | X_1, \ldots, X_k) = \frac{\prod_{i=1}^k p(X_i | \theta)}{\int_{\Gamma_j} \prod_{i=1}^k p(X_i | \theta) \, d\pi(\theta)},$$

then

$$\exp(-\rho(1 - \rho) D_\rho^{\mathrm{Re}}(q(X) \| p_j(X)))$$

$$= \mathbf{E}_{X_1,\ldots,X_k} \left( \frac{\int_{\Gamma_j} \prod_{i=1}^k p(X_i|\theta)\,d\pi(\theta)}{\pi(\Gamma_j) \prod_{i=1}^k q(X_i)} \right)^\rho$$

$$\times\, \mathbf{E}_{X_{k+1}} \left( \frac{\int_{\Gamma_j} w(\theta|X_1,\ldots,X_k) p(X_{k+1}|\theta)\,d\pi(\theta)}{q(X_{k+1})} \right)^\rho$$

$$\leq \mathbf{E}_{X_1,\ldots,X_k} \left( \frac{1/\pi(\Gamma_j) \int_{\Gamma_j} \prod_{i=1}^k p(X_i|\theta)\,d\pi(\theta)}{\prod_{i=1}^k q(X_i)} \right)^\rho$$

$$\times\, \sup_{p\in \mathrm{co}(\Gamma_j)} \mathbf{E}_{X_{k+1}} \left( \frac{p(X_{k+1})}{q(X_{k+1})} \right)^\rho$$

$$= \mathbf{E}_{X_1,\ldots,X_k} \left( \frac{1/\pi(\Gamma_j) \int_{\Gamma_j} \prod_{i=1}^k p(X_i|\theta)\,d\pi(\theta)}{\prod_{i=1}^k q(X_i)} \right)^\rho$$

$$\times\, \sup_{p\in \mathrm{co}(\Gamma_j)} e^{-\rho(1-\rho) D_\rho^{\mathrm{Re}}(q(X_{k+1})\|p(X_{k+1}))}$$

$$\leq e^{-\rho(1-\rho)k \inf_{p\in \mathrm{co}(\Gamma_j)} D_\rho^{\mathrm{Re}}(q(X_1)\|p(X_1))} \cdot \sup_{p\in \mathrm{co}(\Gamma_j)} e^{-\rho(1-\rho) D_\rho^{\mathrm{Re}}(q(X_{k+1})\|p(X_{k+1}))}$$

$$= \exp\left( -\rho(1-\rho)n \inf_{p\in \mathrm{co}(\Gamma_j)} D_\rho^{\mathrm{Re}}(q(X_1)\|p(X_1)) \right).$$

This proves the claim for $n = k+1$. Note that in the above derivation, the first of the two inequalities follows from the fact that with fixed $X_1,\ldots,X_k$, the density $p(X_{k+1}) = \int_{\Gamma_j} w_i(\theta|X_1,\ldots,X_k) p(X_{k+1}|\theta)\,d\pi(\theta) \in \mathrm{co}(\Gamma_j)$; the second of the two inequalities follows from the induction hypothesis. $\square$

PROOF OF THEOREM 5.3. We simply substitute the estimates of Lemma A.4 into Lemma A.3. $\square$

**Acknowledgments.** The author would like to thank Andrew Barron for helpful discussions that motivated some ideas presented in this paper, and Matthias Seeger for useful comments on an earlier version of the paper. The author would also like to thank the anonymous referees for helpful comments and for pointing out related papers.

## REFERENCES

[1] BARRON, A. and COVER, T. (1991). Minimum complexity density estimation. *IEEE Trans. Inform. Theory* **37** 1034–1054. MR1111806

[2] BARRON, A., SCHERVISH, M. J. and WASSERMAN, L. (1999). The consistency of posterior distributions in nonparametric problems. *Ann. Statist.* **27** 536–561. MR1714718

[3] CATONI, O. (2004). A PAC-Bayesian approach to adaptive classification. Available at www.proba.jussieu.fr/users/catoni/homepage/classif.pdf.

[4] GHOSAL, S., GHOSH, J. K. and VAN DER VAART, A. W. (2000). Convergence rates of posterior distributions. *Ann. Statist.* **28** 500–531. MR1790007

[5] LE CAM, L. (1973). Convergence of estimates under dimensionality restrictions. *Ann. Statist.* **1** 38–53. MR0334381

[6] LE CAM, L. (1986). *Asymptotic Methods in Statistical Decision Theory.* Springer, New York. MR0856411

[7] LI, J. (1999). Estimation of mixture models. Ph.D. dissertation, Dept. Statistics, Yale Univ.

[8] MEIR, R. and ZHANG, T. (2003). Generalization error bounds for Bayesian mixture algorithms. *J. Mach. Learn. Res.* **4** 839–860.

[9] RÉNYI, A. (1961). On measures of entropy and information. *Proc. Fourth Berkeley Symp. Math. Statist. Probab.* **1** 547–561. Univ. California Press, Berkeley. MR0132570

[10] RISSANEN, J. (1989). *Stochastic Complexity in Statistical Inquiry.* World Scientific, Singapore. MR1082556

[11] SEEGER, M. (2002). PAC-Bayesian generalization error bounds for Gaussian process classification. *J. Mach. Learn. Res.* **3** 233–269. MR1971338

[12] SHEN, X. and WASSERMAN, L. (2001). Rates of convergence of posterior distributions. *Ann. Statist.* **29** 687–714. MR1865337

[13] VAN DE GEER, S. (2000). *Empirical Processes in M-Estimation.* Cambridge Univ. Press.

[14] VAN DER VAART, A. W. and WELLNER, J. A. (1996). *Weak Convergence and Empirical Processes. With Applications to Statistics.* Springer, New York. MR1385671

[15] WALKER, S. and HJORT, N. (2001). On Bayesian consistency. *J. R. Stat. Soc. Ser. B Stat. Methodol.* **63** 811–821. MR1872068

[16] YANG, Y. and BARRON, A. (1999). Information-theoretic determination of minimax rates of convergence. *Ann. Statist.* **27** 1564–1599. MR1742500

[17] ZHANG, T. (1999). Theoretical analysis of a class of randomized regularization methods. In *Proc. Twelfth Annual Conference on Computational Learning Theory* 156–163. ACM Press, New York. MR1811611

[18] ZHANG, T. (2004). Learning bounds for a generalized family of Bayesian posterior distributions. In *Advances in Neural Information Processing Systems 16* (S. Thrun, L. K. Saul and B. Schölkopf, eds.) 1149–1156. MIT Press, Cambridge, MA.

YAHOO RESEARCH
135 OAKLAND AVENUE
TUCKAHOE, NEW YORK 10707
USA
E-MAIL: tzhang@yahoo-inc.com